\documentclass[mmnp]{edpsmath}
\usepackage{graphicx}

\usepackage[none]{hyphenat}
\usepackage{comment}
\usepackage{caption}
\usepackage{wrapfig}
\usepackage{exscale}
\usepackage{amsmath}
\usepackage{multicol}
\usepackage{tikz,graphicx,caption,subcaption} 
\usepackage{amsthm}
\graphicspath{{fig/}}
\usepackage{tikz,graphicx,caption,subcaption} 

\tikzset{
	invisible/.style={opacity=0,text opacity=0},
	visible on/.style={alt=#1{}{invisible}},
	alt/.code args={<#1>#2#3}{%
		\alt<#1>{\pgfkeysalso{#2}}{\pgfkeysalso{#3}} % \pgfkeysalso doesn't change the path
	},
}
\usepackage{float}
\usepackage{booktabs}
\usetikzlibrary{snakes}
\usetikzlibrary{decorations.pathmorphing}
\usepackage[shortlabels]{enumitem}
\usepackage{graphicx,caption}
\usepackage{textcomp}
\usepackage{adjustbox}
\usepackage{siunitx}
\usepackage{amssymb}
\newcommand{\BibTeX}{{\rm B\kern-.05em{\sc i\kern-.025em b}\kern-.08em\ignorespaces
  T\kern-.1667em\lower.7ex\hbox{E}\kern-.125emX}}

\begin{document}
\overfullrule=10pt
%%%%%%%%%%%%%%%%%%%%%%%%%%%%%%%%%%%%%%%%%%%%%%%%%%%%%%%%%%%%%%%%%%
%                         THE TOP MATTER                         %
%%%%%%%%%%%%%%%%%%%%%%%%%%%%%%%%%%%%%%%%%%%%%%%%%%%%%%%%%%%%%%%%%%
%
\title{Mathematical modeling of leukemia chemotherapy in bone marrow}%
\thanks{This work has been partially supported by the Fundaci\'on Española para la Ciencia y la Tecnolog\'ia (FECYT project PR214), the Asociaci\'on Pablo Ugarte (APU, Spain), Junta de Andaluc\'ia (Spain) group FQM-201, Ministry of Science and Technology, Spain (grant number PID2019-110895RB-I00, funded by MCIN/AEI/ 10.13039/501100011033). This work was also subsidized by a grant for the research and biomedical innovation in the health sciences within the framework of the Integrated Territorial Initiative (ITI) for the province of Cadiz (grant number ITI-0038-2019). ITI is $80\%$ co-financed by the funds of the FEDER Operational Program of Andalusia 2014-2020 (Council of Health and Families).}

\runningtitle{The \textsl{edpsmath} documentclass users guide Version 2 }

\author{Ana Niño-López}
\address{Department of Mathematics, Universidad de C\'{a}diz, Puerto Real, C\'{a}diz, Spain;
\email{ana.nino@uca.es}}
\secondaddress{Biomedical Research and Innovation Institute of Cadiz (INiBICA) Hospital Universitario Puerta del Mar, C\'adiz, Spain}
\author{Salvador Chulián}
\sameaddress{1,2}
\author{Álvaro Martínez-Rubio}
\sameaddress{1,2}
\author{Cristina Blázquez-Goñi}
\sameaddress{2}\secondaddress{Department of Pediatric Hematology and Oncology, Hospital de Jerez, C\'adiz, Spain}
\author{María Rosa}
\sameaddress{1,2}
\runningauthors{The \textsl{edpsmath} documentclass users guide Version 2}
\begin{abstract}

Acute Lymphoblastic Leukemia (ALL) accounts for the $80\%$ of leukemias when coming down to pediatric ages. Survival of these patients has increased by a considerable amount in recent years. However, around $15-20\%$ of treatments are unsuccessful. For this reason, it is definitely required to come up with new strategies to study and select which patients are at higher risk of relapse.  Thus the importance to monitor the amount of leukemic cells to predict relapses in the first treatment phase.

\vspace{2ex}
In this work we develop a mathematical model describing the behavior of ALL, examining the evolution of a leukemic clone when treatment is applied. In the study of this model it can be observed how the risk of relapse is connected with the response in the first treatment phase. This model is able to simulate cell dynamics without treatment, representing a virtual patient bone marrow behavior. 
Furthermore, several parameters are related to treatment dynamics, therefore proposing a basis for future works regarding childhood ALL survival improvement.

\end{abstract}
\begin{resume}
La leucémie aiguë lymphoblastique (LAL) représente $80\%$ des leucémies
diagnostiquées chez l'enfant. La survie de ces patients a considérablement
augmenté ces dernières années. Malgré cela, environ $15$ à $20\%$ des traitements échouent. Pour cette raison, il est crucial de réaliser de nouvelles stratégies d’étude
et de sélectionner les patients présentant le risque de rechute le plus élevé. D’où
l’importance de contrôler le nombre de cellules leucémiques pour prédire les
rechutes dans la première phase du traitement.\newline

Dans ce travail, un modèle mathématique est développé qui décrit le comportement
de la LAL, en étudiant l’évolution du clone leucémique lorsque le traitement est
appliqué. Dans l’étude de ce modèle, nous pouvons observer comment le risque de
rechute est lié à la réponse à la première phase du traitement. Ce modèle est
capable de simuler la dynamique cellulaire sans traitement, représentant le
comportement de la moelle osseuse d’un patient virtuel. De plus, les différents
paramètres liés à la dynamique du traitement sont étudiés. Par conséquent, la base
des travaux futurs visant à améliorer la survie de la LAL chez l´enfant est proposée.

\end{resume}
\maketitle

%%%%%%%%%%%%%%%%%%%%%%%%%%%%%%%%%%%%%%%%%%%%%%%%%%%%%%%%%%%%%%%%%%
%                              THE TEXT                          %
%%%%%%%%%%%%%%%%%%%%%%%%%%%%%%%%%%%%%%%%%%%%%%%%%%%%%%%%%%%%%%%%%%
\tableofcontents
\section*{Introduction}
Leukemia is a malignant disease originating in the bone marrow. 
Particularly, it arises from a disruption in hematopoiesis, the process in charge of blood cells production \cite{orkin2008hematopoiesis,jagannathan2013hematopoiesis}.
Hematopoiesis is usually depicted as hierarchical tree, in which a hematopoietic stem cell (HSC) can differentiate into other cells from the lymphoid or myeloid line. Depending on the branch, platelet, red blood cells and lymphocytes are produced. Leukemia is not only distinguished by the linage but also depending on the maturation stage of the cancer cells. Acute Lymphoblastic Leukemia (ALL) is caused by cells with fast growth in the lymphoid branch, with special incidence in pediatric patients \cite{pui2001childhood, pui2008acute}.

\vspace{2ex}

Each type and subtype of leukemia has an associated protocol that specifies the therapeutic recommendations, and that varies from country to country \cite{egler2016asparaginase}. In particular, SEHOP-PETHEMA-2013 protocol is used in Spain to treat pediatric patients with ALL. This protocol was designed by Spanish Society of Pediatric Hematology and Oncology along with Spanish Hematology Treatment Program and it is regularly reviewed \cite{mesegue2021lower,ruiz2022venous}.

There are different treatments and medications depending on the risk group assigned to a patient: standard, intermediate or high. In such protocols, patients go through successive phases: induction, consolidation, re-induction and maintenance. Each phase consists of the combination of chemotherapeutic agents and corticoids with varying schedules. As reflected in other works \cite{hunger2012improved, ma2014survival}, the progressive improvement and modification of these protocols has increased survival rates to around $80$ to $85\%$. Nonetheless, $15$ to $20\%$ of patients fail to achieve long term remission and therefore relapse \cite{ward2019estimating}. Apart from adding to and modifying these protocols, recent reviews have highlighted the need for new therapies and approaches that would tackle the disease differently \cite{bhojwani2013relapsed,terwilliger2017acute}. This also includes the improvement of risk assignment, the management of secondary effects and the study of personalized medicine, whose power lies not only in treatment but also in prevention \cite{mathur2017personalized}. \\

In this sense, mathematical models have recently found a wide variety of applications in medicine \cite{bocharov2018mathematical,tosenberger2013modelling}. In particular, the discipline of mathematical oncology has been able to approach recent oncological issues \cite{gatenby2003mathematical,altrock2015mathematics}. To understand cell dynamics, several works have described biological processes related to cellular behaviour \cite{mackey1994global, komarova2013principles} and also related to hematopoiesis and lymphopoiesis \cite{bessonov2006cell,lochem2004, stiehl2011characterization}. This has allowed the mathematical characterization of the hierarchical structure of blood cell lineages  \cite{bonnet1997human, anderson2011genetic, marciniak2009,salvi2020} and their development in bone marrow. 
Specifically, ALL has been considered in previous mathematical models \cite{marciniak2019, clapp2015review,ducrot2007model} with an structural vision of leukemic cells similar to that of healthy cells. In order to fulfill their potential, these models should be validated against appropriate treatment data \cite{moricke2005prognostic,dai2021clinical,bullhallmarks2022}. 
Based on the current protocols to treat ALL patients  \cite{pui2006treatment,ronghe2001remission}, some mathematical models describe the drugs' effect regardless of the leukemia type \cite{jayachandran2014optimal, mouser2014model},  while other focus explicitly on the myeloid linage \cite{nave2022new, rubinow1976mathematical,pefani2014chemotherapy}. Other mathematical models describe alternative strategies to chemotherapy, such as CAR T-cell therapy \cite{cartcells,perez2021car, leon2021car,kimmel2021roles}, which is successful in patients that do not respond to the standard treatment.\\

Inspired by the examples above, we set out in this study to build a mathematical model that describes the evolution of ALL and its response to treatment in the bone marrow. This type of leukemia has received less attention from the mathematical community, and can benefit just as much. The mathematical model would allow the derivation of theoretical results (like stability analysis) but also clinical information like the relative importance of parameters, the estimation of response times or the influence of dosage and schedule. It could eventually be used to fit individual patient data and obtain insights into the possibility of personalized therapy for this disease \cite{egler2016asparaginase,clapp2015review}.\\

In this paper, we focus on the simulation of standard risk ALL pediatric patients in the first phase of the treatment, or induction phase, which lasts 28 days. In terms of drugs, this phase comprises the corticosteroid Prednisone, the chemotherapeutic agents Vincristine and Daunorubicine, and Asparaginase, each with their own dosages and infusion times. In addition, triple intrathecal is given to prevent central nervous system disease involvement. The modelling of this treatment stage requires an analysis of the bone marrow environment in addition to the assumptions related to the growth of leukemic cells. We study this scenario and simulate how treatment affects both leukemia and the lymphoid linage in bone marrow. Besides, we propose a classification of the patients depending on their leukemic cells development to predict which patients will respond to treatment correctly.\\

The work has been organised as follows: In Sec. \ref{sec:mat&met} we describe an existing healthy lymphopoiesis mathematical model and extend it to include leukemia. This extended model will also include treatment administration. In Sec. \ref{sec:results} we show model simulations of leukemia growth and its evolution under therapy, obtaining realistic ranges for parameters and highlighting the most relevant ones in terms of response. In Sec. \ref{sec:discussion} we discuss the results obtained. %and analyze the importance of progressing in this field.

\section{Material and methods}
\label{sec:mat&met}
\begin{figure}[!ht]
  \centering
 % \begin{adjustbox}{max width=\textwidth,center}
      \includegraphics[width=0.9\textwidth]{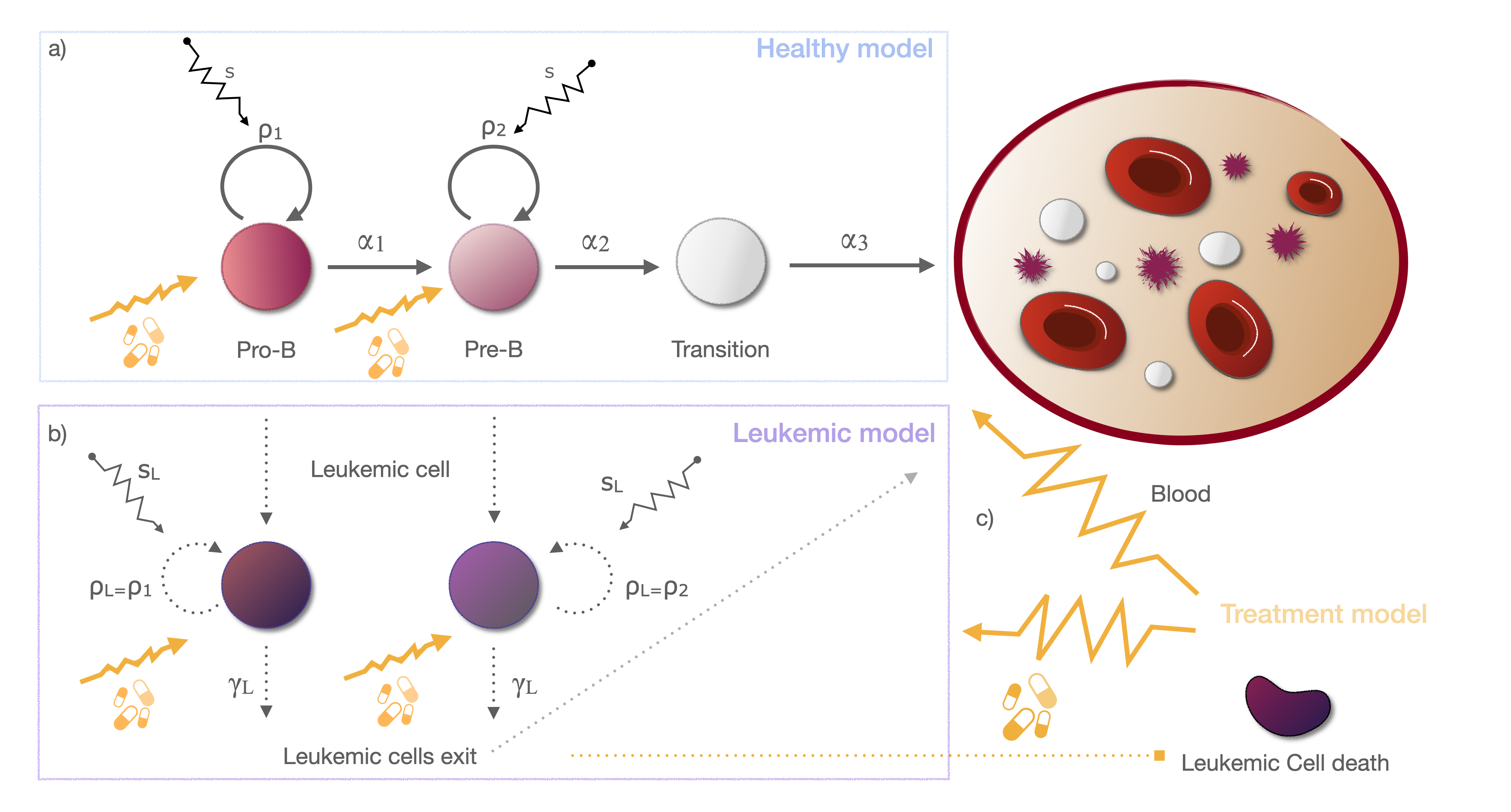}
   %   \end{adjustbox}
     \caption{\textbf{Diagram of cell development along with the appearance of leukemia and its treatment.} \textbf{a) } Immature healthy B cells from the bone marrow can develop into more mature cells as follows: In an early maturation stage, a Pro-B cell can either proliferate in its own compartment with a rate $\rho_1$ or progress to the next compartment (Pre-B) with rate $\alpha_1$. Analogously, a Pre-B cell can remain in a Pre-B stage with a proliferation rate $\rho_2$ and develop to the transition compartment with rate $\alpha_2$. Finally we consider that transition cells do not proliferate in their own compartment and they arrive to the blood flux with rate $\alpha_3$. Growth rate is regulated by the feedback signal $s$.   \textbf{b)} Leukemic cells could come from a Pro-B stage or Pre-B stage, inheriting the proliferation rate $\rho_L$ of the related compartment. A regulatory signal $s_L$ controls leukemic cell development. Furthermore, for the case of ALL we assume that leukemic cells do not differentiate into a next stage and therefore $\gamma_L$ indicates the blood exit rate for leukemic cells into blood. \textbf{c) } When drug is administered, both healthy and leukemic cells are affected by the treatment, with cells dying at a rate proportional to their proliferation rate.}
    \label{fig:figure1}
  \end{figure}

In this section we describe a mathematical model of healthy, homeostatic lymphopoiesis in the bone marrow. We then include the appearance of a leukemic cell and study its evolution. We finally approximate a treatment function following the current treatment scheme during induction phase and model its effect in the previous models. An overview of all such models is shown in Figure~\ref{fig:figure1}.

\subsection{Healthy bone marrow model}
\label{subsection:salvi}
Following the immunophenotypical characterization of B cells \cite{lochem2004}, previous mathematical models, have considered B lymphopoiesis as a sequential process dividing their maturation stages in Pro-B (early stage characterized by immunophenotypic markers CD10\textsuperscript{+}/CD45\textsuperscript{-}), Pre-B (intermediate stage with CD10\textsuperscript{+}/CD45\textsuperscript{+}) and Transition cells (last stage with CD10\textsuperscript{-}/CD45\textsuperscript{+}) \cite{marciniak2009,marciniak2019}. Mathematically, this can be regarded as three compartments respectively depending on each stage: $C_1(t)$, $C_2(t)$ and $C_3(t)$, along with their associated parameters related to proliferation and differentiation. This was translated in \cite{salvi2020} into the compartmental model presented below:
\begin{subequations}
\label{eq:salvi0}
\begin{align}
\dfrac{dC_1}{dt}&= \bar{s}\rho_1C_1-\alpha_1C_1,\\
\dfrac{dC_2}{dt}&= \bar{s}\rho_2C_2+\alpha_1C_1-\alpha_2C_2,\\
\dfrac{dC_3}{dt}&= \alpha_2C_2-\alpha_3C_3,
\end{align}
\end{subequations}

\noindent where $\rho_i$, for $i=1,2$, is the proliferation rate and $\alpha_i$, $i=1,2,3$, is the transition rate related to each compartment $i$.\\

In this model the authors also took into account a signal in charge of regulating cell dynamics, $\bar{s}=\bar{s}(t)$, with the following expression:

\begin{equation}
\label{eq:senal0}
\bar{s}=\bar{s}(t)=\dfrac{1}{1+k N},
\end{equation}

\noindent where $k$ is the inhibitory parameter and $N$ the cells subpopulation in charge of the feedback signalling. According to the results in \cite{salvi2020}, we consider here $N=\sum_{i=1}^3 C_i$.
This signalling was shown to affect predominately the proliferation rate due to biological reasons such as bounded growth or mathematical ones such as positivity of the solutions. Having reviewed this model, we expand it by including the appearance of a leukemic clone.

\subsection{Modeling leukemia cells}
\label{subsec:leuk}
Now, in addition to the three compartments $C_1,C_2,C_3$, we consider a leukemic cells compartment, $L$. We assume that leukemic cells grow as a logistic function where the curve's maximum value is $L_{\max}$ and $\rho_L$ is the logistic growth rate. Just like cells in a healthy hematopoiesis, a leukemic clone is characterized by a maturation stage \cite{bonnet1997human,marciniak2019}. Therefore, we assume $L$ to depend  on the stage where the leukemic clone originates, i.e., from the Pro-B stage ($\rho_L=\rho_1$) or from Pre-B stage ($\rho_L=\rho_2$). Finally, $\gamma_L$ would denote the blood exit rate of leukemic cells. Here we do not consider clones coming from the transition stage due to the fact that ALL is produced by immature, proliferating cells.\\

 We consider then the following equation to define the flux in the leukemic cells compartment:
\begin{equation}
\label{leuk0}
    \dfrac{dL}{dt}=s_L\rho_L L \left(1-\dfrac{L}{L_{\max}}\right) -\gamma_L L.
\end{equation}

Finally, signalling $s_L$ in this compartment has a similar behavior to the healthy bone marrow compartments as in \cite{salvi2020}, and we assume it to affect the proliferation rate as in Model \eqref{eq:salvi0}. Therefore, 
\begin{equation}
\label{eq:sl}
s_L=\bar{s}=\dfrac{1}{1+k \sum_{i=1}^3 C_i},
\end{equation}

\noindent assuming leukemic cells should not be affected by their own signalling due to the evasion of growth suppressors \cite{hallmarks2011}.
Nevertheless, healthy stages are influenced by cells in compartment $L$, due to the fact that the more cells in $L$, the less healthy cells develop due to the invasion of the bone marrow. Consequently, healthy signalling $\bar{s}$ in Model~\eqref{eq:salvi0} should be modified to: 
\begin{equation}
\label{eq:signalsanasleuk}
 s=\dfrac{1}{1+k\left(L+\sum_{i=1}^3{C_i}\right)}.
\end{equation}

In addition to these assumptions, we can also include a constant influx of healthy stem cells $c_0$ in our initial compartment $C_1$ \cite{cartcells}, which results in the following model, including both leukemic and healthy B cells:

\begin{subequations}
\label{eq:ModelosClon}
\begin{align}
\dfrac{dC_1}{dt}&= c_0 + s\rho_1C_1-\alpha_1C_1,\\
\dfrac{dC_2}{dt}&= s\rho_2C_2+\alpha_1C_1-\alpha_2C_2,\\
\dfrac{dC_3}{dt}&= \alpha_2C_2-\alpha_3C_3,\\
\dfrac{d L}{d t} &= s_{L} \rho_{L} L\left(1-\dfrac{L}{L_{\max }}\right)-\gamma_{L} L,
\end{align}
\end{subequations}
with $s$ as in Eq.~\eqref{eq:signalsanasleuk} and $s_L$ as Eq.~\eqref{eq:sl}, where $\rho_L$ depends on the origin of the leukemic clone \cite{anderson2011genetic, marciniak2019}. \\

This model allows us to study different treatment regimes, to compare with actual results and to include and test what is known about the behaviour and effect of these drugs \cite{egler2016asparaginase,ronghe2001remission, karon1966vincris}.

\subsection{Modeling treatment for leukemia}
\label{subsec:treatment}

\begin{figure}[!ht]
  \centering
   \begin{adjustbox}{max width=0.9\textwidth,center}
      \includegraphics[width=\textwidth]{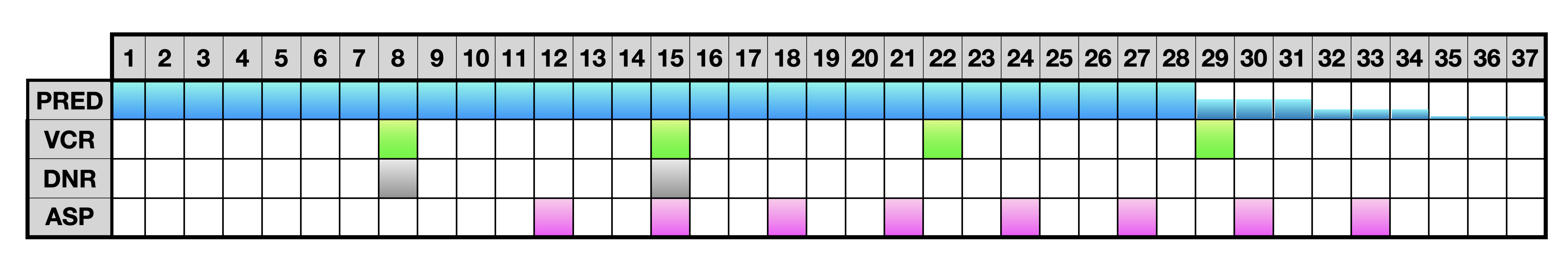}
      \end{adjustbox}
\caption{\textbf{Schedule for Induction I'A treatment for a Standard Risk Patient.} Data from SEHOP-PETHEMA protocol: Prednisone (PRED) is administered at $60mg/m^2/day$ for the first $28$ days. Then, the dose is reduced to $30mg/m^2/day$ for days $29,30$ and $31$, to $15mg/m^2/day$ for days $32,33$ and $34$ and to $7.5mg/m^2/day$ for days $35,36$ and $37$. Vincristine (VCR): $1.5 mg/day$ on days $8,15,22$ and $29$. Daunorubicine (DNR): $30mg/m^2/day$ on days $8$ and $15$. Asparaginase (ASP): $10000 U/m^2/day$ on days $12,15,18,21,24,27,30$ and $33$.}
     %Total dose of drugs administered to the patient during Induction I'A phase.
    \label{fig:figure2}    
\end{figure}

As explained previously, there are different treatment protocols that can be applied \cite{clapp2015review}, depending on the drugs used, their doses, or even the time of infusion. We consider the Induction I'A phase, as described by the SEHOP- PETHEMA-2013 protocol, for a standard risk patient. This phase consists of four drugs and lasts 37 days. Triple intrathecal is not included in the study since it does not affect the lymphoid branch. Each drug has a schedule with corresponding doses and days of administration, depending on the patient's body surface, as indicated in Figure~\ref{fig:figure2}.

 To simplify, we consider treatment days as a scale which corresponds leukemia detection day with day $+0$ and treatment beginning with day $+1$. The most important days for monitoring the patient are treatment days $+8$, $+15$ and $+33$.  In those days, bone marrow sample extractions are carried out and lymphocyte levels are studied. In day $+8$ a blood extraction is done to assess response to prednisone. According to SEHOP-PETHEMA-2013 protocol, there must be less than $10^6$ leukemic cells per milliliters of blood in day $+8$. The Minimal Residual Disease (MRD) is studied in bone marrow, whose positivity is defined by the presence of $0.01\%$ or more leukemic cells in bone marrow \cite{campana2010minimal}. Therefore, patient responds to the treatment if leukemic cells involve at most $0.01\%$ in bone marrow  in day $+15$, and there are no blasts in day $+33$. \\
 
Following studies on the behaviour of each drug \cite{ drugbankPrednisona, drugbankVincri, drugbankDauno, drugbankAspa, armstrong1987applications}, we consider the treatment $\mu_j(t)=\mu_j$ as a function that describes the amount of drug $j$ in the bone marrow at time $t$,
being $j\in J=\{ P,V,D,A\}$: Prednisone $(\mu_P)$, Vincristine $(\mu_V)$, Daunorubicin $(\mu_D)$ and Asparaginase $(\mu_A)$. Once the dose is administered, the drug has an exponential decrease associated with that medicine half-life.

We then define $\mu_j:\mathbb{R\textsuperscript{+}}\rightarrow \mathbb{R\textsuperscript{+}}$ as

\begin{equation}
\label{eqdifmu}
\dfrac{d\mu_j}{dt}= -\lambda_j \mu_j,\\
\end{equation}

\noindent with $\lambda_j$ related to each drug $j$ half-life $\tau_j$ measured in days:
\begin{equation}\lambda_j= \dfrac{\log(2)}{\tau_j}.\end{equation}
\vspace{1em}
Each drug $j$ has a different dose $q_j$ in several days $\mathcal{D}_j\subset \mathbb{N}$, as shown in Figure \ref{fig:figure2}. We define the dose administered as: 

\begin{equation}
 Q_j(t)=
 \begin{cases} 
      q_j& t \in \mathcal{D}_j,\\
      &\\
      0 & t \notin\mathcal{D}_j.\\
    \end{cases}
\end{equation}

\begin{figure}[!ht]
     \hspace{-3em}\includegraphics[width=0.9\textwidth]{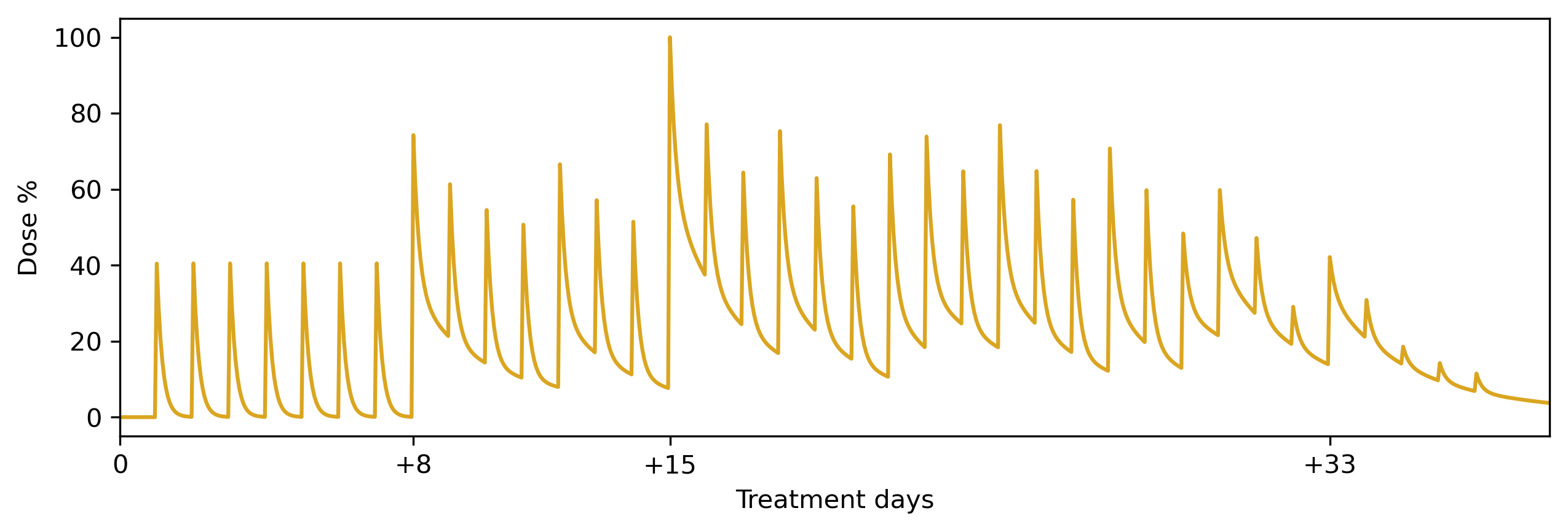}
\caption{\textbf{Total treatment influence in bone marrow}. 
Percentage of total drug dose for Induction phase, following SEHOP-PETHEMA-2013 (Figure \ref{fig:figure2}) and with parameters from Table~\ref{table:parametros}.}
     %Total dose of drugs administered to the patient during Induction I'A phase.
    \label{fig:figure3}    
\end{figure}

We finally consider the total treatment function $\mu$, represented in Figure~\ref{fig:figure3}, as the weighted sum of all drugs:

\begin{equation}
\label{eqmu}
    \mu=\sum_{j\in J}\delta_j(\mu_j+Q_j),
\end{equation}
where $\delta_j$ indicates the influence of drug $j$ on the total effect of the treatment.\\

\begin{comment}

let $q_j(t)=q_j$ be the treatment dose $j$ applied in time $t$ according to SEHOP-PETHEMA-2013 protocol. Here, in the moment of induction, each drug $j$ will be related to a control parameter $\delta_j$ which indicates how each treatment $j$ affects healthy and leukemic cells death. Besides,
$\delta_j,q_j \in \mathbb{R}$, $\delta_j, q_j>0$ and $k_j$ is defined by the half-life of each drug \cite{drugbankPrednisona, drugbankVincri, drugbankDauno, drugbankAspa}. We finally consider the total drug function $\mu$ as the addition of all treatments, with

\begin{equation}
\label{eqmufin}
      \dfrac{d\mu}{dt}=\sum^4_{j=1}\dfrac{d\mu_j}{dt}.
\end{equation}

\end{comment}

Finally, we include a term in Model \eqref{eq:ModelosClon}
that indicates cellular death by chemotherapy. We consider this does not affect transition cells, owing to the fact that they do not have any associated proliferation. The resulting model is shown below ~\eqref{eq:ModeloTratamiento}.

%Consequently, we need to include into the Model \eqref{eq:ModelosClon} a term which indicates the cellular death by treatment specifically $-\mu\rho_iC_i$ for $i=1,2$. We consider no cellular death due to treatment for transition cells, owing to the fact that they do not have any associated proliferation. \\

%Finally, the Model~\eqref{eq:ModeloTratamiento} below sums up all considerations stated, with

\begin{subequations}
\label{eq:ModeloTratamiento}
\begin{align}
\dfrac{dC_1}{dt}&= c_0 + s\rho_1C_1-\alpha_1C_1-\mu\rho_1C_1,\\
\dfrac{dC_2}{dt}&= s\rho_2C_2+\alpha_1C_1-\alpha_2C_2-\mu\rho_2C_2,\\
\dfrac{dC_3}{dt}&= \alpha_2C_2-\alpha_3C_3,\\
\dfrac{d L}{d t}&= s_{L} \rho_{L} L\left(1-\dfrac{L}{L_{\max }}\right)-\gamma_{L} L-\mu \rho_L L,\\
%las ecuaciones del tto
\mu&= \sum_{j\in J}\delta_j \left( \mu_j + Q_j\right),\\
Q_j&=\begin{cases} 
      q_j& t \in \mathcal{D}_j,\\
      &\\
      0 & t \notin\mathcal{D}_j,\\
   \end{cases},\\
 \dfrac{d\mu_j}{dt}&=
      -\lambda_j \mu_j(t).
\end{align}
\end{subequations}

\vspace{1em}
Particularly, when $\mu=0$, (i.e. no treatment is applied) Model~\eqref{eq:ModeloTratamiento} transforms into Model~\eqref{eq:ModelosClon}.

\subsection{Parameters estimation}
%From the data found in the literature, we adjusted the parameters for our assumptions. Thus, we work from the data collected in the Table ~\ref{table:parametros}.
  We consider literature data for healthy bone marrow properties \cite{salvi2020} along with parameters related to leukemia \cite{cartcells}. Drugs half lives are taken from pharmacokinetics and pharmacodynamics studies. All parameter values are included in Table~\ref{table:parametros}. The estimations for dosage assume a $25-30$ $kg$ body mass for a child, which in turn implies a body surface area equivalent to $1 \; m^2$ \cite{haycock1978geometric}.\\

\begin{table}[!ht] 
%%% \tablesize{} %% You can specify the fontsize here, e.g., \tablesize{\footnotesize}. If commented out \small will be used.
\begin{adjustbox}{max width=1.1\textwidth,center}
\begin{tabular}{ccccc}
\toprule
\textbf{Parameter}	& \textbf{Meaning}	&
%\textbf{Model 1 Value} & 
\textbf{Value} & \textbf{Unit}  & \textbf{Source} \\
\midrule
$c_0$& Influx cells from HSC & $10^7$ & $cell$ & Estimated from \cite{cartcells}\\
$\rho_1$ & Pro-B proliferation rate &
%$0.6931$&
$\log (2)$&$day^{-1}$& \cite{salvi2020}\\
 
 $\rho_2$ & Pre-B proliferation rate &
 %$0.4621$&
 $\log(2)/1.5$&$day^{-1}$& \cite{salvi2020}\\

$\alpha_1$ & Transition rate: Pro-B to Pre-B %&$0.168$
& $0.168$&$day^{-1}$& \cite{salvi2020}\\

$\alpha_2$ & Transition rate: Pre-B to Transition  &
%$0.144$&
$0.144$&$day^{-1}$& \cite{salvi2020}\\

$\alpha_3$ & Blood exit rate &
%$0.288$&
$0.288$&$day^{-1}$& \cite{salvi2020}\\
 $k$ &Signal Intensity&
 %$10^{-10}$&
 $10^{-10}$&$cell^{-1}$&  \cite{cartcells}\\
%\hline
%$k_L?$ &$10^{-10}$&$10^{-10}$\\
%\hline
$\gamma_L$ & Leukemic cells blood exit rate & %$0.000288$&
$0.288 \times 10^{-3}$&$day^{-1}$& \cite{cartcells}\\

$L_{\max}$ &Leukemic cells carrying capacity& %$10^{12}$&
$10^{12}$&$cell$& \cite{cartcells}\\
\midrule
$q_P$ & Prednisone dose &
$60$& $mg/day$ & {\small SEHOP-PETHEMA-2013}\\

$q_V$ & Vincristine dose & 
$1.5$&$mg/day$ & {\small SEHOP-PETHEMA-2013}\\

$q_D$ & Daunorubicin dose & 
$30$&$mg/day$ & {\small SEHOP-PETHEMA-2013}\\

$q_A$ & Asparaginase dose & 
$10000$&$U/day$ & {\small SEHOP-PETHEMA-2013}\\
&&&&\\
$\delta_P$ & Prednisone influence & %$10^{12}$&
%$0.06$
$\left[\dfrac{1}{60},\dfrac{1}{6}\right]$&  $day/mg$&Estimated  \\
&&&&\\
$\delta_V$ & Vincristine influence & %$10^{12}$&
%$0.93$
$\left[\dfrac{1}{1.5},\dfrac{10}{1.5}\right]$&  $day/mg$&Estimated \\

&&&&\\
$\delta_D$ & Daunorubicin influence & 
%$0.1$
$\left[\dfrac{1}{30},\dfrac{1}{3}\right]$& $day/mg$ &Estimated   \\

&&&&\\
$\delta_A$ & Asparaginase influence & %$10^{12}$&
%$2.5 \times 10^{-4}$
$\left[\dfrac{1}{10^{4}},\dfrac{1}{10^3}\right]$&  $day/U$&Estimated \\
&&&&\\

$\lambda_P$ & Prednisone decaying rate & 
$9.6 \log(2)$&$day^{-1}$  & \cite{drugbankPrednisona}\\
$\lambda_V$ & Vincristine decaying rate & 
$0.28 \log(2)$&$day^{-1}$  & \cite{drugbankVincri}\\
$\lambda_D$ & Daunorubicin decaying rate & 
$1.17 \log(2)$&$day^{-1}$  & \cite{drugbankDauno}\\
$\lambda_A$ & Asparaginase decaying rate & 
$0.8\log(2)$&$day^{-1}$ & \cite{drugbankAspa}\\

\midrule
$C_1(0)$ & Pro-B cells & 
%$7.1565\times 10^9$&
$3.52211\times 10^9$&$cell$& \cite{salvi2020}\\
 $C_2(0)$ & Pre-B cells & 
 %$8.70833\times 10^{10}$&
 $1.84911\times10^{10}$&$cell$& \cite{salvi2020}\\
$C_3(0)$ & Transition cells &  %$2.6125\times10^{10}$&
$9.24555\times10^9$&$cell$& \cite{salvi2020}\\
$L(0)$ & Leukemic cells & 
%$1$&
$1$&$cell$& Assumption \\

\bottomrule
\end{tabular}
\end{adjustbox}
\caption{\textbf{Parameter values}. $L_{\max}$ is measured in number of cells and $k$ in $cell^{-1}$ while the rest of the presented parameters related to the proliferation and transition rates are measured in $day^{-1}$. On the other hand, medicine dosage $q_j$ is measured in $mg/day$ and the influence of each drug is considered $day/mg$. Asparaginase is a particular case since is measured in $\mathcal{U}/day$ and its influence in $day/\mathcal{U}$. Values of $\lambda_j$ are measured in $day^{-1}$ due to the fact that they come from half-life values.  Values used to initialize each case $C_i(0)$, $i=1,2,3$, have been obtained from the steady values in healthy bone marrow models \cite{salvi2020} along with our parameters. A single leukemic cell is supposed at the beginning, from which Model~\eqref{eq:ModelosClon} starts, and the influx of cells $c_0$ is also obtained from Model~\eqref{eq:salvi0} steady states.}
\label{table:parametros}
\end{table}

The study of all estimated parameters is presented in Supplementary Information document. Control parameters related to drug, $\delta_j$ with $j\in J$, have been selected by searching for values that implies that the patient responds to treatment correctly, i.e. considering an amount of blasts less than $10^6$ $blasts/ml$ in blood in day $+8$ and MRD$<0.01\%$ in bone marrow in day $+15$ ~\cite{szczepanski2007mdr}. In preliminary simulations we find that if $\delta_jq_j <10^0$, the dose administered is ineffective. We further consider that $\delta_jq_j>10^1$ involves the death of the patient due to drug toxicity. We then set all $\delta_j q_j$ to the same order of magnitude $\left[10^0, 10^1\right]$, $\forall j \in J$. In this way we ensure Model~\eqref{eq:ModeloTratamiento} can represent a virtual patient whose behavior is similar to real patients according to the amount of cells in each bone marrow aspiration.

\section{Results}
\label{sec:results}
\subsection{Leukemic and healthy cells dynamics without treatment are properly represented by the model.} \label{resultado1:leukemia}

\begin{figure}[!ht]
  \centering
  
     \hspace{-2em} \includegraphics[width=\textwidth]{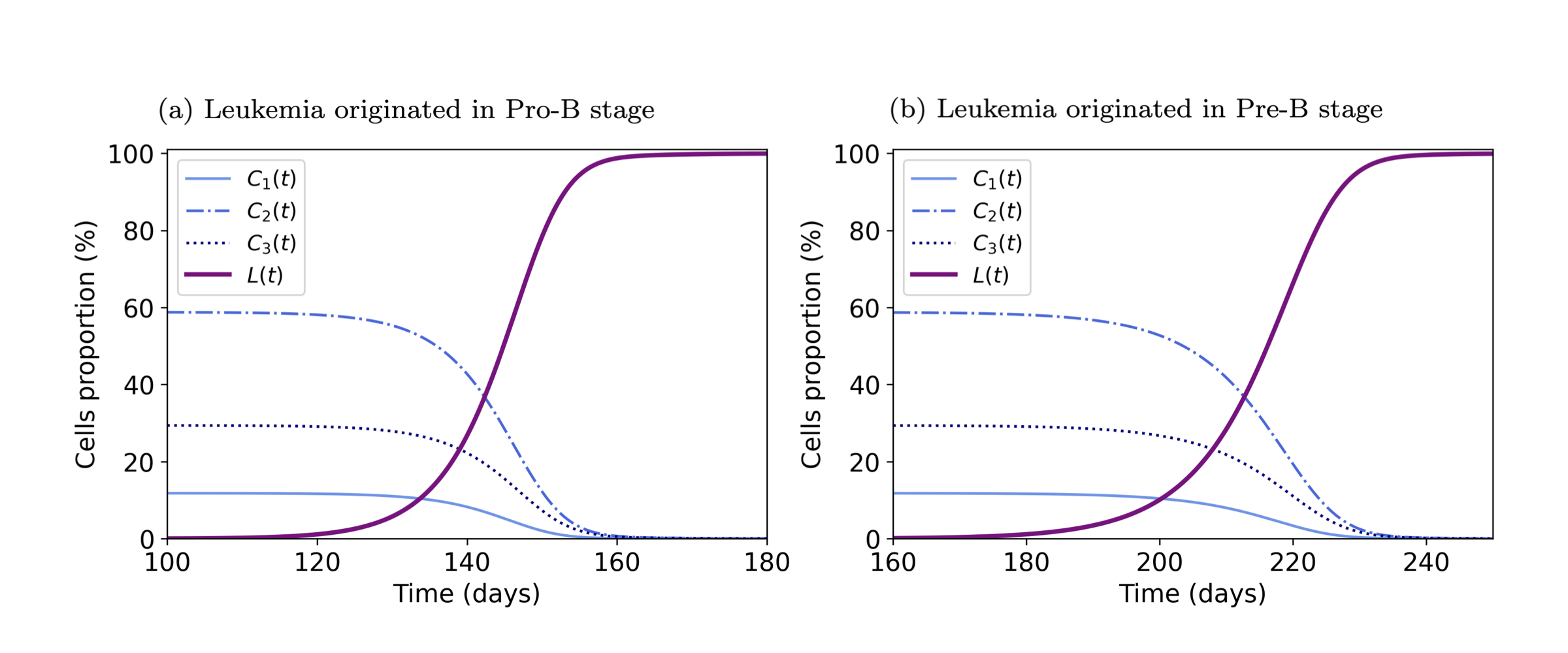}

     \caption{\textbf{Dynamics of cell compartments proportion in the bone marrow in the presence of the leukemic clone.} Initial data and parameters from Table~\ref{table:parametros}. (\textbf{a}) Case $\rho_L=\rho_1=0.6931\ day^{-1}$. Leukemic cells, $L$, (solid purple line) have a logistic growth and reach $80\%$ around day $150$. At the same time, Pro-B cells ($C_1$, solid blue line), Pre-B cells ($C_2$, dash-dotted blue line) and Transition cells ($C_3$, dotted blue line) decrease. (\textbf{b}) Case $\rho_L=\rho_2=0.4621\ day^{-1}$. Dynamics is similar to the first case but leukemia is detected around the day 223.}
    \label{fig:figure4}
  \end{figure}

%The first point it must be taken into account is healthy bone marrow stability values which are initial values in this model \cite{salvi2020}.
We now consider no treatment is applied, Model~\eqref{eq:ModelosClon}, i.e., Model ~\eqref{eq:ModeloTratamiento} with $\mu_j=0$, $\forall j \in J$. This model combines the healthy dynamics of immature B-cells and how leukemic cells stability values as initial inputs for the model. We simulate the appearance of a leukemic cell (at $t=0$), depending on the leukemic clone originates, either in a Pro-B $\left(\rho_L=\rho_1\right)$ or a Pre-B stage $\left(\rho_L=\rho_2\right)$. Results of solving in a range to $\left[0,300\right]$ days are shown in Figure~\ref{fig:figure4}, displaying the range of days with more information in each case.\\

We next show the simulations for each case. Figure~\ref{fig:figure4} represents Model~\eqref{eq:ModelosClon} cells evolution when $\rho_L=\rho_1$ is considered for the first case, and secondly, model for $\rho_L=\rho_2$ is simulated. In addition, for each time $t$ measured in days, the bone marrow cells proportions are shown for the different maturation compartments. 

\begin{comment}

  \begin{subfigure}[b]{0.48\textwidth}
      \caption{Leukemia originated in Pro-B stage}
      \includegraphics[width=\textwidth]{fig/LeukemiaModelP1.png}
    \end{subfigure}
        \begin{subfigure}[b]{0.48\textwidth}
      \caption{Leukemia originated in Pre-B stage} \includegraphics[width=\textwidth]{fig/LeukemiaModelP2.png}
    \end{subfigure}

\end{comment}

As it can be observed in Figure~\ref{fig:figure4}, both graphics represent the total bone marrow invasion by leukemia. Depending on the stage from which the first leukemic cell comes, leukemic cells proliferate more or less quickly. The appearance of a leukemic clone is set to be detected when the leukemic population accounts for $80\%$ of blasts from the total population of B lymphocytes, \cite{amin2005having}. When the first leukemic cell is originated in Pro-B stage, $\rho_L=\rho_1$, the disease is detected around day $150$. However, if the original leukemic cell comes from the Pre-B stage, leukemia is appreciated from day $223$. Both two cases imply a progressive decreasing in healthy cells due to the fact that the invasive ability of leukemic cells prevents the normal development of cells in the bone marrow.

\subsection{The first blood extraction in day $+8$ of treatment allows to approximate the influence of prednisone ($\delta_P$).}
\label{sec:prednisoneinfluence}
%When leukemia is detected, there are approximately $80\%$ of leukemic cells in bone marrow since data collected from bone marrow samples that day of different patients confirm it \cite{amin2005having}. It occurs around 150 days since the first leukemic cell appeared. Therefore, 
Treatment begins according to treatment schedule in SEHOP-PETHEMA-2013 protocol, after $80\%$ of blasts are detected, producing a drop in both leukemic and healthy cells.
%It occurs when there are approximately $80\%$ of leukemic cells in bone marrow since data collected from bone marrow samples that day of different patients confirm it. \\

In order to find the optimal parameters to represent Model~\eqref{eq:ModeloTratamiento}, we pay attention to the first blood extraction (day $+8$). To do so, we estimate relations between blasts in bone marrow and blood, owing to the fact that there is a correlation between them \cite{bianconi2013estimation,choi2014reference}.\\

\begin{figure}[!ht]
 \includegraphics[width=1\textwidth]{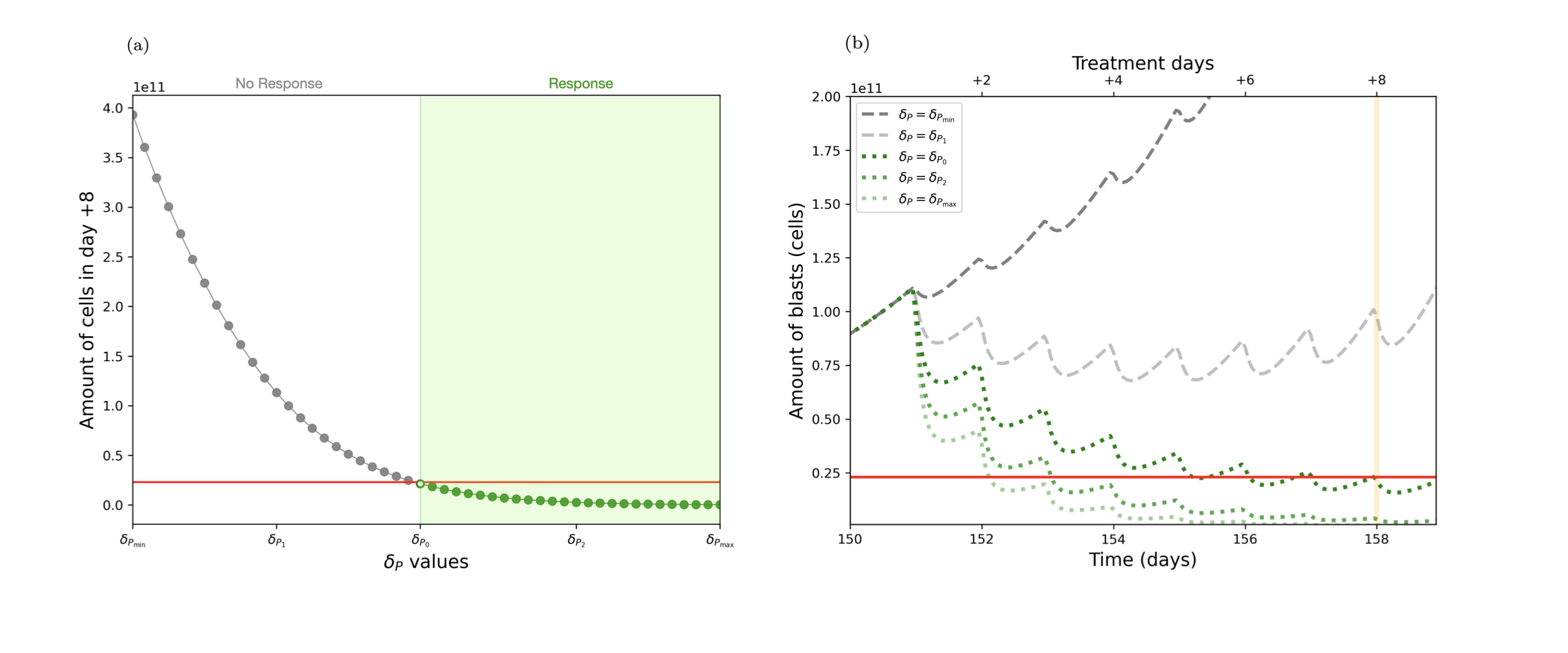}
   \caption{\textbf{Variation of $\delta_P$ values to analyze results in day $+8$ of treatment.} As presented in Table \ref{table:parametros}, $\delta_P\in\left[1/60,1/6\right]=\left[\delta_{P_{\min}},\delta_{P_{\max}}\right]$. 
   Five values of $\delta_P$ are highlighted: the midpoint between the maximum and minimum values $\delta_{P_{0}}\approx0.09$, $\delta_{P_{1}}\approx 0.05$ the middle between $\delta_{P_{\min}}$ and $\delta_{P_{0}}$, and $\delta_{P_{2}}\approx 0.13$ the center of the interval $\left[\delta_{P_{0}},\delta_{P_{\max}}\right]$ .   (a) Model \eqref{eq:ModeloTratamiento} is solved for $50$ values of prednisone influence. For each $\delta_P$ value, amount of blasts in day +8 is obtained and associated with a patient who responses (green points) or not (gray points) to prednisone. (b) Leukemic cells evolution depending on several $\delta_P$ values. Green-dotted lines denote patients who response to prednisone, while gray-dashed lines represent bone marrow dynamics which do not reach any response to prednisone from data in day $+8$ (yellow vertical line). Horizontal solid red line represents the limit of blasts assumed by SEHOP-PETHEMA-2013 protocol to consider a good response of less than $2.3\times10^{10}$ blasts in bone marrow.}
  \label{fig:figure5}
\end{figure}

SEHOP-PETHEMA-2013 protocol defines a good response to prednisone if there is fewer than $10^6$ blasts per milliliter of blood. We assume a pediatric patient with $25-30 kg$ body mass, and due to the blood - body mass correlation  \cite{raes2006reference, linderkamp1977estimation}, this assumption entails approximately $2.1-2.5$ liters of blood. We consider a patient with $2.3l$ of blood to our study, and consequently, patients should have less than $2.3\times 10^9$ blasts in blood. It implies that leukemic cells limit in the bone marrow should be less than $2.3\times10^{10}$ blasts at day $+8$ to obtain a good response to prednisone. We analyze all $\delta_P$ possible values within our range in Table \ref{table:parametros}, taking into account leukemic cells evolution for the first eight days of treatment (see Figure \ref{fig:figure5}).\\

As it is shown in Figure \ref{fig:figure5} (a), the amount of blasts in bone marrow in day $+8$ of treatment implies a direct consequence related to possible $\delta_P$ values. 
Values lower than $\delta_{P_0}$ imply that patients do not response to prednisone, while values higher than $\delta_{P_0}$ entail a good prognosis for the patient due to the fact that in day $+8$ of treatment reach less than the limit of blasts assumed.\\

On the other hand, in Figure \ref{fig:figure5} (b) the development of leukemic cells according to Model~\eqref{eq:ModeloTratamiento} is depicted for five different values of $\delta_P$. Values which allow to reach blasts levels lower than $2.3\times10^{10}$ total blasts in day $+8$ corresponds to the green area in Figure \ref{fig:figure5} (a). 

%\textbf{Como vemos en la Fig a para que haya una buena respuesta del tratamiento y que funcione, debe considerarse un deltaP a partir de aproximadamente 0.09. como vemos en los ejemplos mostrados en (b) para los valores menores a ese delta p no se tiene una buena respuesta del paciente y para valores mayores sí. Con todo esto, la tabla de parámetros se puede afinar y reajustar los valores  de deltaP que cumplan esta condición.  }

%En los artículos \cite{raes2006reference, linderkamp1977estimation} da la correlación entre la cantidad de sangre y el peso en niños. El segundo es del 77.\\

%Se tiene por tanto que para un niño de 25-30 kg se tienen entre 2100 y 2500 ml de sangre. Supongamos 2300ml de sangre para los cálculos. En el día +8 según pethema deberá haber menos de $10^6 blastos/ml$ por lo tanto debe haber menos de $2.3\times 10^9$ blastos totales en sangre. Relacionando sangre y médula como que en médula hay de un orden de magnitud más de células\cite{bianconi2013estimation,choi2014reference}, se tiene por tanto que el número de blastos totales en médula no puede ser mayor que  $2.3\times 10^{10}$ para poder afirmar que el paciente responde correctamente a la prednisona.

\subsection{Prednisone and Vincristine influence values analysis implies a classification of patients}

Let us how consider all parameters $\delta_j$ in order to analyze their behaviour in Model \eqref{eq:ModeloTratamiento}, taking into account the number of blasts in bone marrow at day $+15$.\\

Based on the parameter sensitivity analysis (Supplementary Information), main influential parameter is $\delta_P$. To study their influence, we vary parameters $\delta_P$ and $\delta_V$, related, respectively, to prednisone and vincristine, and check how many leukemic cells are there day $+15$ of treatment schedule in comparison to the healthy population (MRD). We fix parameters $\delta_D=2.5/30 day/mg$ and $\delta_A=2.5/10^4 day/U$ given their minimal influence in the sensitivity analysis. 
We focus on the variation of $\delta_P$ values in comparison with $\delta_V$ values to obtain the minimal value which implies the patient responds to treatment. We assume patient responds to treatment correctly if  day $+15$ of treatment leukemic cells represent less than $0.01\%$ of the whole lymphocyte population.\\
%after 100 days of their administration in comparison to the healthy population. 
%We focus on the variation of $\delta_P$ values in comparison with $\delta_V$ values to obtain the minimal value which implies no relapse. We assume no relapse if 100 days after treatment's end, leukemic cells represent less than $0.01\%$ of the whole lymphocyte population.\\
%According to our assumptions, if the leukemic cells percentage in bone marrow is fewer than $0.01\%$ 100 days after the treatment finished, there is no relapse.\\

\begin{figure}[!ht]
   \centering
      \includegraphics[width=0.7\textwidth]{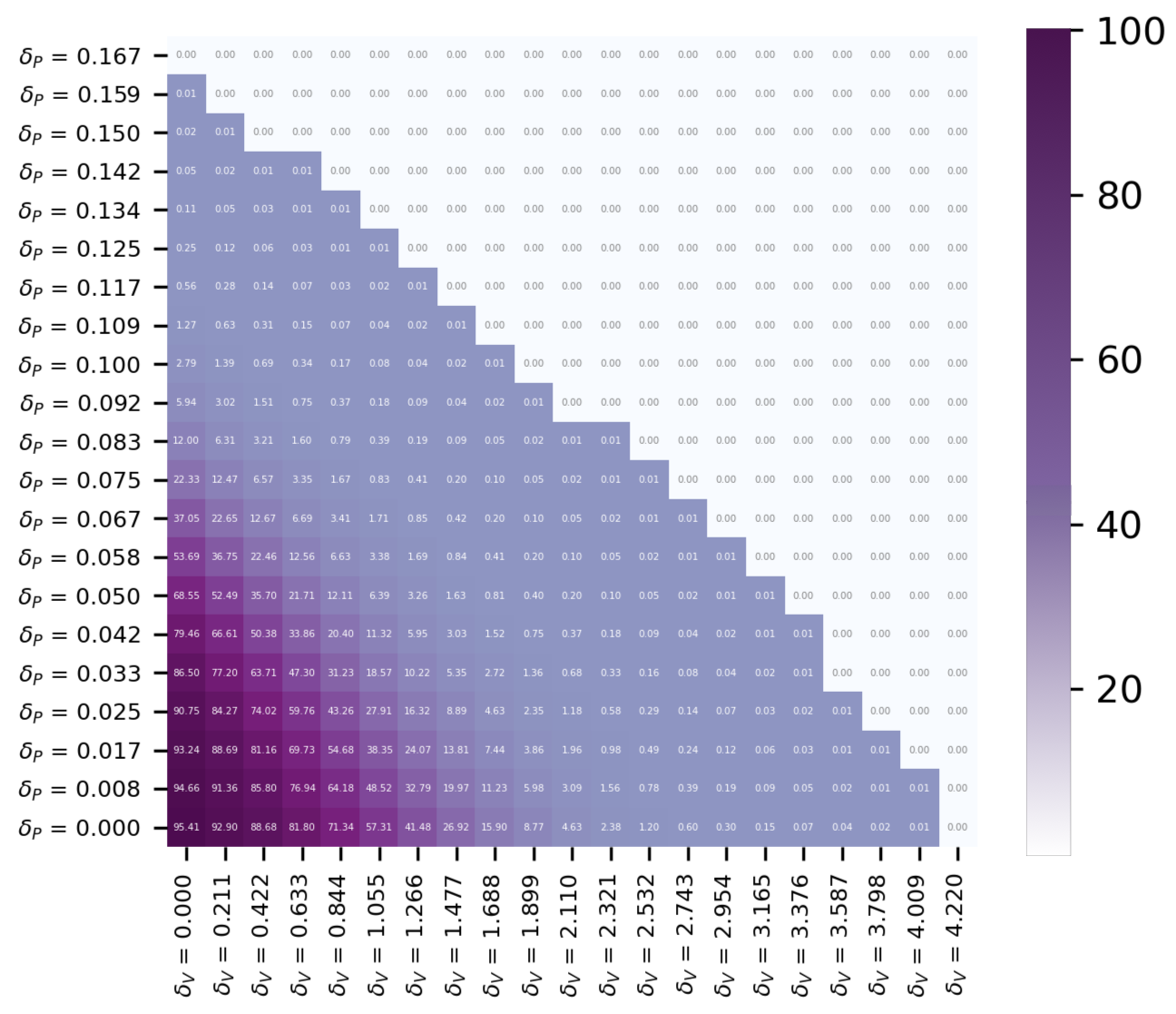}
     \caption{\textbf{Blasts percentage in bone marrow day $+15$ of treatment.} We consider $21$ values up to $0.167 day/mg$ for Prednisone and $21$ values up to $4.22 day/mg$ for Vincristine. For each values combination, Model \eqref{eq:ModeloTratamiento} is solved for parameters in Table \ref{table:parametros} with $\delta_D= 2.5/30 day/mg$, $\delta_A= 2.5/10^4 day/U$. The value in each box corresponds to the percentage of leukemic cells in day $+15$ depending on values $\delta_P$ and $\delta_V$.}
    \label{fig:figure6}
  \end{figure}

Values ranges considered for $\delta_j$, $j\in J$ are those in Table \ref{table:parametros}. In Supplementary Information, we show a heat map which takes parameters in the widest range. That range of parameter for such heat map has been adapted to a more meaningful range. 
We show in Figure~\ref{fig:figure6} a heat map which connects parameters variation with the blasts percentage in day $+15$.

%\subsection{Study of drug influence allows a classification of patients}
This heat map leads to associate all values with two semiplanes which will be able to predict if the patient responds correctly depending on the values which match his data. 
Regarding Figure ~\ref{fig:figure6}, we can approximate the influence of each drug in a patient and predict if the patient responds to the treatment correctly or not. Firstly, parameters related to the influence of drugs in cells death, $\delta_j$, can be approximated for a set of data from a specific patient in several days of evolution of the treatment. Thus, we classify these parameters according to Figure~\ref{fig:figure6} in purple or white zone to predict which percentage of leukemic cells will be there in bone marrow day $+15$. There will be a good prognosis if leukemic compartment entails less than $0.01\%$ (MDR) in bone marrow, white zone, while if those values are in the purple zone, the patient will not respond correctly to treatment. \\

According to choose the minimal values which implies a good prognosis of the patient, $\delta_P=0.092day/mg$ and $\delta_V=2.11day/mg$ will be parameters used in the next results.
 
\begin{comment}

This can be illustrated briefly by Figure~\ref{fig:variacionP}, in which three cases are represented. As it is observed there, prednisone influence affects results after 100 days, obtaining a relapse or a recovery depending on that value. 

   \begin{figure}[!ht]
  \centering
          \includegraphics[width=0.6\textwidth]{fig/variacionDeltaPfeb.png}
     \caption{\textbf{Variation of leukemia dynamics along with therapy according to $\delta_P$ values.} Leukemic cells have different behavior depending on the considered values of prednisone influence. Model solved with data from Table~\ref{table:parametros} $\delta_V= 1.15$ $day/mg$, $\delta_D=2.5/30$ $day/mg$, $\delta_A=2.5\times10^{-4}$ $day/U$, obtaining the resulting graphics: blue, orange and green lines, depending on $\delta_P$ values, $0.01, 0.04$  or $0.09$  $day/mg$, respectively. }
    \label{fig:variacionP}
  \end{figure}

\end{comment}

\subsection{A new model about bone marrow behavior in presence of treatment is proposed to study.}

%From the leukemic cell appearance simulations, we could simulate some types of treatment as it is said before. In this paper, we consider treatment from SEHOP-PETHEMA-2013 shown in Figure ~\ref{fig:figure2}. 

\begin{figure}[!ht]
  \centering
   
    \hspace{-2em}\includegraphics[width=1\textwidth]{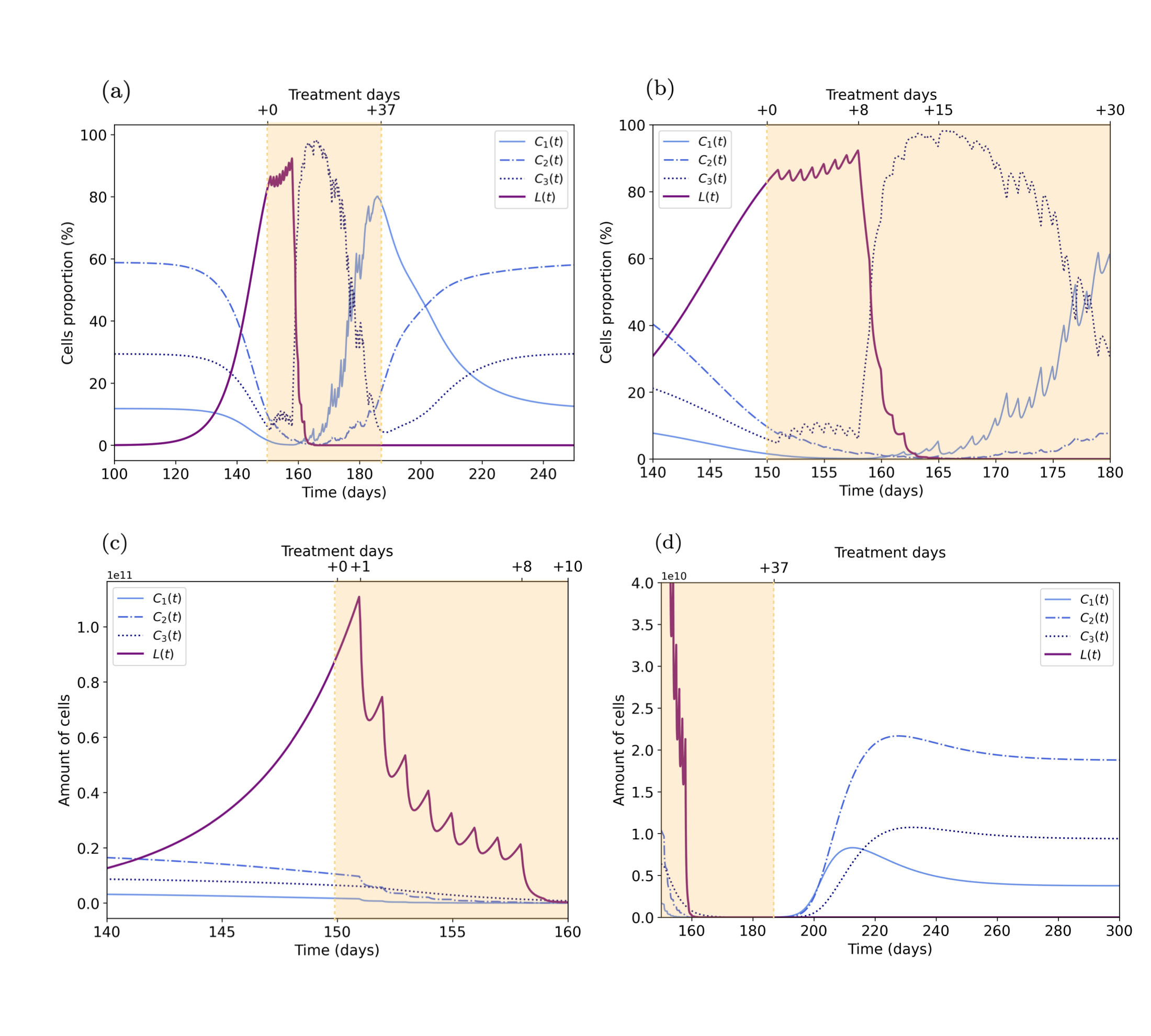}

     \caption{\textbf{Dynamics of cell compartments in the bone marrow in the presence of the leukemic clone along with therapy.} Initial data, parameters related to leukemia and parameters about treatment are those presented in Table~\ref{table:parametros} except for $\delta_P=0.092day/mg$, $\delta_V=2.11day/mg$, $\delta_D=2.5/30 day/mg$, $\delta_A=2.5/10^4 day/U$ and $\rho_L=\rho_1=0.6931\ day^{-1}$ case which we consider without loss of generality. \textbf{(a)} Cells proportion evolution in bone marrow from the first leukemic cell appearance along with the treatment application. \textbf{(b)} Bone marrow dynamics centered between day $150$ and day $180$, corresponding to $+0$ and $+30$ from therapy, respectively. \textbf{(c)} Compartment values measured in amount of cells in ten first days of treatment. It is observed that day $+8$ is the most aggressive of them. \textbf{(d)} Treatment finishes in day $+37$. The extinction of leukemic cells implies healthy cells back to normal.}
    \label{fig:figure7}
  \end{figure}

After proposing Model~\eqref{eq:ModeloTratamiento}, we simulate the situation using data in Table~\ref{table:parametros} along with parameters associated to the $\delta_j$ values study in previous sections. As proposed in Sec. \ref{subsec:leuk} it is possible to observe that leukemic cells have a higher proliferation than the rest of lymphocytes. Moreover, as said before, it is known that the blasts level in bone marrow when leukemia is detected is around $80\%$. We assume that treatment starts the next day and present results in Figure~\ref{fig:figure7}.\\

%When leukemia is detected, there are approximately $80\%$ of leukemic cells in bone marrow since data collected from bone marrow samples that day of different patients confirm it \cite{amin2005having}. It occurs around 150 days since the first leukemic cell appeared. Therefore, treatment begins producing a drop in both leukemic and healthy cells. It occurs when there are approximately $80\%$ of leukemic cells in bone marrow since data collected from bone marrow samples that day of different patients confirm it.\\

%On the 15th day after treatment beginning (day $+15$), leukemic cells proportion are adjusted to the samples, obtaining values lower than MRD. Furthermore, on the day $+33$, it is very important do not appreciate leukemic cells to proceed with the transplant.\\

Therefore, according to the available data about different patients, we have been able to approximate the evolution of bone marrow cells in the presence of treatment as it is shown in Figure ~\ref{fig:figure7}. It can be observed that the treatment affects both healthy and leukemic cells and that once the disease is eliminated, healthy cells compartments return to normal levels.
\clearpage

Generally, in Figure~\ref{fig:figure7}, a leukemic bone marrow with treatment is represented. From usual cells compartments proportions, leukemia develops until leukemic cells exceed $80\%$ of total occupancy in bone marrow. From those conditions, the therapy is applied for $37$ days, and then, bone marrow normal levels are recovered. In a reduced range of days (Figure~\ref{fig:figure7}(b)) it can be appreciated how cells proportions are modified according to drugs administered each day. In fact, day $+8$ is the first aggressive day in which leukemic cells do not grow so fast, therefore, in following days this population continues its reduction. It is for that reason that the highest cells proportion in bone marrow for about $30$ days of treatment is the transition compartment, $C_3$, which are not influenced by treatment since they do not proliferate. Paying attention to the amount of cells (Figure~\ref{fig:figure7}(c)-(d)), it is possible to verify that all of the cells in bone marrow decrease gradually and finally, return to normal levels. 

%Firstly, we consider leukemia to be detected when leukemic cells constitute approximately $80\%$ of the total population of the B lymphocytes. Without loss of generality, we consider the case $\rho_L=\rho_1$. Henceforth, treatment is applied (yellow shaded area), altering each compartment value, depending on the influence regarding the proliferation rate as it is presented in Model ~\eqref{eq:ModeloTratamiento}. It can be appreciated how it is reached the MRD in ALL \cite{campana2010minimal,szczepanski2007mdr} after 15 days, as indicated in SEHOP-PETHEMA-2013.

\section{Discussion}
\label{sec:discussion}
Leukemia is the cancer with the highest incidence in pediatric age, particularly, Acute Lymphoblastic Leukemia and $20\%$ of these diagnoses do not respond correctly to treatment, and as a consequence, there is a relapse. This problem is the main reason of our study and the principal question to answer is why the treatment does not work for those patients. Our intention is to find patterns in bone marrow along with medicines applied, so that, it is possible to predict relapses and improve current treatments. We approximate treatment influence in cell death to analyze the behavior of each bone marrow based on flow cytometry data throughout treatment time.
Modeling current treatment in leukemia B leads to advances towards optimizing the application of medicines. After studying leukemic cells behavior it is possible to analyze and advance in this area. Available clinical data are important to determinate parameters values according to different bone marrow samples from the same patient. \\

First of all, we have reviewed a healthy bone marrow model in Sec.~\ref{subsection:salvi}. That study allows us to include the assumption of the appearance of a leukemic cell. It is found that leukemic cells have a strong development due to their properties of malignant disease. The leukemic cell develops as it has been described in Model \eqref{eq:ModelosClon}, and it could be more or less accelerated depending on the inherited proliferation rate for $\rho_L$, as shown in Figure~\ref{fig:figure4}. In fact, available data \cite{moricke2005prognostic, amin2005having} of hematologic patients reinforce that the disease is diagnosed when there are about $80\%$ leukemic cells in bone marrow. As it is mentioned in Sec.~\ref{resultado1:leukemia}, six steady states are obtained by studying Model~\eqref{eq:ModelosClon}. The unique stable steady state is the one that represents the total occupancy of the bone marrow by leukemic cells. The rest of solutions are unstable. There are several possible explanations for this result.  On one hand, both two solutions which have negative compartments, are biologically senseless results. In spite of the fact that there are two more solutions which are non negative, they lack biological meaning. Moreover, there is a solution which reproduces healthy bone marrow behavior without leukemia, Sec.~\ref{subsection:salvi}. Therefore, Model~\eqref{eq:ModelosClon} replicates reality and its study allows to adjust parameters related to leukemic population growth. 
In the case in which leukemia appears in the Pro-B stage,  $\rho_L=\rho_1$ for Model~\eqref{eq:ModelosClon}, the amount of leukemic cells is the same as Pro-B cells around 127 days; it is equal to Transition cells in the day 132; and Pre-B cells are reached by leukemic cells in 135 days.
Simulations for model which assume the appearance of the leukemic cell in Pre-B stage, $\rho_L=\rho_2$, have a similar behavior, but the process is slower.
In general, both cases have the same behavior. Every simulation in the presence of the leukemic clone allow for deducing that healthy and leukemic cells coexist for a period of time which depends on the case, and then, leukemic cells take up all the bone marrow capacity, therefore, there is no space for healthy cells.
According to the Model \eqref{eq:ModelosClon} and its simulations, it is important to take into account that cells development is supposed without any treatment which is responsible for reducing the leukemic cells proportion. Hence, it explains why there are that suddenly changes three months later. \\

From the information collected in the SEHOP-PETHEMA-2013 protocol, we model the treatment behavior in the bone marrow when it is administered. Without loss of generality, we consider a leukemic cell originated in Pro-B compartment. 
From the first leukemic cell appearance, the proliferation of that population takes about 150 days to be diagnosed with the result that when it is detected, the disease reproduces more quickly. That is why treatment begins to be administered in the following days, thus, there are no data related to the same patient in different times without treatment. We have focused on a standard risk patient in
Induction I'A which is the first phase of the treatment and it is a determining phase in the therapy. If there are signs which do not correspond to proper treatment response, the protocol for this patient is modified and associated to other risk level.
%Therefore, we consider a virtual patient who responds properly to medication in its first phase. Prednisone, Vincristine, Daunorubicine and Asparaginase levels are expressed by Eq.~\eqref{eqmu} which involve cells death in bone marrow throught Model~\eqref{eq:ModeloTratamiento}.
Therefore, we consider a therapy based on Prednisone, Vincristine, Daunorubicine and Asparaginase, administered according to the schedule shown in Figure~\ref{fig:figure2}. Taking into account the total amount of drugs in the body during the treatment, we assume an exponential decay of each drug from its application in order to different references about timing and dosage drugs effects \cite{clapp2015review,acharya1984development,kay2013improving}. In addition, we include the term $\delta_j$ as the influence of each drug in cells death, to obtain the function $\mu$ which imply leukemic and healthy cells death by treatment depending on their proliferation. 

When a sample from a patient is studied as a possible diagnosis of ALL is owing to the patient has previously developed some of the hallmarks of cancer \cite{hallmarks2011}. Generally, bone marrow first extractions indicate that leukemia takes up around $80\%$ of bone marrow. Therefore, protocol is immediately activated. That treatment approximation implies an improvement about the perception related to drugs actions in the body. Depending on each patient, the influence of each drug will be more or less aggressive, hence, the importance of the personalized medicine \cite{mathur2017personalized}.  It is for that reason that a detailed study about drugs influences must be contemplated to study those behaviors. The aim is to provide each new patient with the best treatment. A patient is considered with good prognosis if some conditions are reached. On the treatment day $+8$, blood sample should show less than $10^6$ blasts per milliliters of blood, on the day $+15$, leukemic cells must entail values smaller than $0.01\%$ and the day $+33$ blasts in bone marrow must not be observed, considering that this is the moment of the transplant. \\

A theoretical model of the bone marrow dynamics, Model~\eqref{eq:ModeloTratamiento}, is presented to study, based on previous works along with new assumptions about leukemic cells behaviour and treatment.
In order to study parameters related to prednisone, we have focused on the first eight days of treatment. In those days only prednisone is administered, therefore, it is possible to approximate optimal $\delta_P$ value which implies the good prognosis condition in day $+8$. From range of $\delta_P$ possible values in Table~\ref{table:parametros}, Model~\eqref{eq:ModeloTratamiento} is solved for 50 values in that range, obtaining amount of blasts in blood in day $+8$ related to each one. It is obtained a value $\delta_{P_0}\approx 0.09$ from which, values imply blasts levels smaller than the limit of amount of cells, consequently, we obtain a range of $\delta_P$ values in which patient response to treatment. This result imply two important aspects to consider. On one hand, a patient whose sample is extracted in day $+8$ could be associated to the corresponding $\delta_P$ and predict the evolution of the disease. On the other hand, if we analyze blood extraction in other day, $+5$ for instance, we could be able to associated $\delta_P$ value before day $+8$ and predict results in that day. Thus, it would be possible to change treatment protocol on time. \\

Apart from the study of the day $+8$ conditions, we analyze good response conditions in day $+15$.
According to sensitivity study, prednisone and vincristine are the most significant parameters and hence, we propose a study about influence parameters related to them. On day $+15$, MRD should be less than $0.01\%$ to obtain good results in patients, therefore, we search for values which lead to those leukemic cells levels when Model~\eqref{eq:ModeloTratamiento} is solved for those values. Leukemic compartment percentages in bone marrow on day $+15$ are obtained for each case and presented in Figure~\ref{fig:figure6} in a heat map. As it is observed, there are two semiplanes corresponding to patients who response and patients who do not response correctly to treatment. With the purpose of minimizing both $\delta_P$ and $\delta_V$, we take average minimal values which imply a good response. What is interesting in this result is the fact that the $\delta_P$ value obtained by this method is the same value found by the study of the blasts levels in day $+8$. This finding reinforces that studies must be based on fitting influence parameters to obtain more predictive information of each patient. \\

\hspace{0.5ex} Finally, simulations supported by previous results in this work are presented. Previously, parameter 
estimation is performed and the bone marrow behavior is studied assuming a leukemic cell appearance. The disease is diagnosed when it entails an important part of cells in bone marrow and it is treated according to protocol. As it is observed in the Figure~\ref{fig:figure7} and generally, treatment application leads to patients recovery and the stabilization of healthy cells levels from fifty days after treatment finish. Treatment is applied when leukemia is diagnosed, that it is supposed in day $150$ from the first leukemic cell appearance since there are about $80\%$ of blasts in bone marrow. During eight days prednisone is applied and it is possible to observe a drop in all leukemic and healthy cells but in day $+8$ is more significant because vincristine and daunorubicin are added to the therapy. From day $+12$ in which asparaginase is administered, leukemic cells values are insignificant in comparison to healthy cells. The fact that this situation occurs agree with parameters choice which guarantee that MRD is less than $0.01\%$ on day $+15$. Graphs interpretation is facilitated by the comparison between amount of cells and proportions. In those simulations, it should be highlighted that leukemic cells do not appear the following $150$ days. For that reason, it is considered a partial patient recovery, not completely due to the logistic growth of leukemia. Moreover, the disease always tend to appear again and that idea explains why treatment should continue with other phases which drop the possibility of any rest of leukemic cells develops.  \\

Being limited to available data, this study lacks accuracy related to drug influence estimated parameters. On one hand, an exhaustive study of each patient would required to be able to model the leukemia growth without treatment. To do so, two bone marrow must be extracted, to observe how cells develop and obtain the parameter related to leukemic proliferation, $\rho_L$, in each patient. 
This option is unfeasible since two bone marrow extractions are rare and unnecessary, as drugs are to be administered to the child as soon as possible. Hence, we propose a future study to relate blood blasts and search for a possible correlation.
Another limitation of our analysis would be the simplicity of the treatment function. Drugs evolution inside the body has a complicated and more detailed process which has been reduced and it produces a lack of information. Absorption time  or action mechanism are some of the conditions to be taken into account in order to get closer to real data. Moreover, only the Induction A phase has been formulated instead of the full treatment. We have focused on this phase of treatment since effects of the next stages (Induction B and Consolidation) are not visible enough. In addition, Induction A is the most interesting phase to study owing to its critical influence throughout treatment. \\

Notwithstanding with these limitations, the study suggests that it is possible to approximate  contemplated parameters related to treatment provided that parameters values for leukemia are corroborated. Additionally, different compartment levels in patients bone marrows with ALL corresponds to values obtained from this previous study in which a virtual patient is simulated.
These findings suggest several courses of action for improving techniques used by protocols. There is a definite need for personalizing medicine. We suggest, therefore, that studies consider mechanisms of action of drugs along with an exhaustive monitoring to prove how they inhibit healthy and leukemic cells proliferation as well as how they destroy cells and the environment. We propose a clinical study about model with those conditions. Comparisons between models including new assumptions will be able to offer some approximations to specific models capable of obtaining realistic parameters. Several patients data are essential to validate the results presented. The effectiveness of that resource will allow to predict relapses in patients who have been diagnosed with ALL in a short period of time after treatment begins. The challenge now is to model full treatment, which takes into account available data and suitable parameters which lead to new research to improve current protocols and, consequently, survival expectancy in childhood ALL.\\

In conclusion, in this paper, a mathematical model of the appearance of leukemic cells in a bone marrow has been presented based on a mathematical model of cells development in a healthy bone marrow. Therefore, we have studied leukemic cells behavior to characterise a model describing the lymphopoiesis process when a cancer cell appears.Furthermore, we have simulated the model proposed and obtained results which coincide with biological basis of cancer. Afterwards, treatment application is included according to current protocol, meaning healthy and leukemic cells death. Dynamics in bone marrow are modeled and studied based on parameters estimation to simulate a virtual patient who responds to treatment.This model provides a basis of future action protocols of treatments of ALL.

\bibliographystyle{elsarticle-num} 
\bibliography{biblio}

\begin{thebibliography}{10}
\expandafter\ifx\csname url\endcsname\relax
  \def\url#1{\texttt{#1}}\fi
\expandafter\ifx\csname urlprefix\endcsname\relax\def\urlprefix{URL }\fi
\expandafter\ifx\csname href\endcsname\relax
  \def\href#1#2{#2} \def\path#1{#1}\fi

\bibitem{orkin2008hematopoiesis}
S.~H. Orkin, L.~I. Zon, Hematopoiesis: an evolving paradigm for stem cell
  biology, Cell 132~(4) (2008) 631--644.

\bibitem{jagannathan2013hematopoiesis}
M.~Jagannathan-Bogdan, L.~I. Zon, Hematopoiesis, Development 140~(12) (2013)
  2463--2467.

\bibitem{pui2001childhood}
C.-H. Pui, D.~Campana, W.~E. Evans, Childhood acute lymphoblastic
  leukaemia--current status and future perspectives, The lancet oncology 2~(10)
  (2001) 597--607.

\bibitem{pui2008acute}
C.-H. Pui, L.~L. Robison, A.~T. Look, Acute lymphoblastic leukaemia, The Lancet
  371~(9617) (2008).

\bibitem{egler2016asparaginase}
R.~A. Egler, S.~P. Ahuja, Y.~Matloub, L-asparaginase in the treatment of
  patients with acute lymphoblastic leukemia, Journal of pharmacology \&
  pharmacotherapeutics 7~(2) (2016) 62.

\bibitem{mesegue2021lower}
M.~Mesegu{\'e}, A.~Alonso-Saladrigues, S.~P{\'e}rez-Jaume, A.~Comes-Escoda,
  J.~L. Dapena, A.~Faura, N.~Conde, A.~Catal{\`a}, A.~Ruiz-Llobet,
  E.~Zapico-Mu{\~n}iz, et~al., Lower incidence of clinical allergy with
  peg-asparaginase upfront versus the sequential use of native e. coli
  asparaginase followed by peg-asp in pediatric patients with acute
  lymphoblastic leukemia, Hematological Oncology 39~(5) (2021) 687--696.

\bibitem{ruiz2022venous}
A.~Ruiz-Llobet, S.~Gassiot, E.~Sarrate, J.~Zubicaray, J.~L. Dapena, S.~Rives,
  J.~Sevilla, {\'A}.~Men{\'a}rguez~L{\'o}pez, M.~Panesso~Romero, C.~Montoya,
  et~al., Venous thromboembolism in pediatric patients with acute lymphoblastic
  leukemia under chemotherapy treatment. risk factors and usefulness of
  thromboprophylaxis. results of lal-sehop-pethema-2013, Journal of Thrombosis
  and Haemostasis (2022).

\bibitem{hunger2012improved}
S.~P. Hunger, X.~Lu, M.~Devidas, B.~M. Camitta, P.~S. Gaynon, N.~J. Winick,
  G.~H. Reaman, W.~L. Carroll, Improved survival for children and adolescents
  with acute lymphoblastic leukemia between 1990 and 2005: a report from the
  children's oncology group, Journal of clinical oncology 30~(14) (2012) 1663.

\bibitem{ma2014survival}
H.~Ma, H.~Sun, X.~Sun, Survival improvement by decade of patients aged 0--14
  years with acute lymphoblastic leukemia: a seer analysis, Scientific reports
  4~(1) (2014) 1--7.

\bibitem{ward2019estimating}
Z.~J. Ward, J.~M. Yeh, N.~Bhakta, A.~L. Frazier, R.~Atun, Estimating the total
  incidence of global childhood cancer: a simulation-based analysis, The Lancet
  Oncology 20~(4) (2019) 483--493.

\bibitem{bhojwani2013relapsed}
D.~Bhojwani, C.-H. Pui, Relapsed childhood acute lymphoblastic leukaemia, The
  Lancet Oncology 14~(6) (2013) e205--e217.

\bibitem{terwilliger2017acute}
T.~Terwilliger, M.~Abdul-Hay, Acute lymphoblastic leukemia: a comprehensive
  review and 2017 update, Blood Cancer Journal 7~(6) (2017) e577--e577.

\bibitem{mathur2017personalized}
S.~Mathur, J.~Sutton, Personalized medicine could transform healthcare,
  Biomedical reports 7~(1) (2017) 3--5.

\bibitem{bocharov2018mathematical}
G.~Bocharov, V.~Volpert, B.~Ludewig, A.~Meyerhans, et~al., Mathematical
  immunology of virus infections, Vol. 245, Springer, 2018.

\bibitem{tosenberger2013modelling}
A.~Tosenberger, F.~Ataullakhanov, N.~Bessonov, M.~Panteleev, A.~Tokarev,
  V.~Volpert, Modelling of thrombus growth in flow with a dpd-pde method,
  Journal of theoretical biology 337 (2013) 30--41.

\bibitem{gatenby2003mathematical}
R.~A. Gatenby, P.~K. Maini, Mathematical oncology: cancer summed up, Nature
  421~(6921) (2003) 321--321.

\bibitem{altrock2015mathematics}
P.~M. Altrock, L.~L. Liu, F.~Michor, The mathematics of cancer: integrating
  quantitative models, Nature Reviews Cancer 15~(12) (2015) 730.

\bibitem{mackey1994global}
M.~C. Mackey, R.~Rudnicki, Global stability in a delayed partial differential
  equation describing cellular replication, Journal of Mathematical Biology
  33~(1) (1994) 89--109.

\bibitem{komarova2013principles}
N.~L. Komarova, Principles of regulation of self-renewing cell lineages, PloS
  one 8~(9) (2013) e72847.

\bibitem{bessonov2006cell}
N.~Bessonov, L.~Pujo-Menjouet, V.~Volpert, Cell modelling of hematopoiesis,
  Mathematical Modelling of Natural Phenomena 1~(2) (2006) 81--103.

\bibitem{lochem2004}
E.~G. Van~Lochem, V.~H.~J. Van~der Velden, H.~K. Wind, J.~G. Te~Marvelde,
  N.~A.~C. Westerdaal, J.~J.~M. Van~Dongen, Immunophenotypic differentiation
  patterns of normal hematopoiesis in human bone marrow: Reference patterns for
  age-related changes and disease-induced shifts, Cytometry Part B: Clinical
  Cytometry 60~(1) (2004) 1--13.

\bibitem{stiehl2011characterization}
T.~Stiehl, A.~Marciniak-Czochra, Characterization of stem cells using
  mathematical models of multistage cell lineages, Mathematical and Computer
  Modelling 53~(7-8) (2011) 1505--1517.

\bibitem{bonnet1997human}
D.~Bonnet, J.~E. Dick, Human acute myeloid leukemia is organized as a hierarchy
  that originates from a primitive hematopoietic cell, Nature medicine 3~(7)
  (1997) 730--737.

\bibitem{anderson2011genetic}
K.~Anderson, C.~Lutz, F.~W. Van~Delft, C.~M. Bateman, Y.~Guo, S.~M. Colman,
  H.~Kempski, A.~V. Moorman, I.~Titley, J.~Swansbury, et~al., Genetic
  variegation of clonal architecture and propagating cells in leukaemia, Nature
  469~(7330) (2011) 356--361.

\bibitem{marciniak2009}
A.~Marciniak-Czochra, T.~Stiehl, A.~D. Ho, W.~J{\"a}ger, W.~Wagner, Modeling of
  asymmetric cell division in hematopoietic stem cells regulation of
  self-renewal is essential for efficient repopulation, Stem cells and
  development 18~(3) (2009) 377--386.

\bibitem{salvi2020}
S.~Chuli{\'a}n, A.~Mart{\'\i}nez-Rubio, A.~Marciniak-Czochra, T.~Stiehl, C.~B.
  Go{\~n}i, J.~F.~R. Guti{\'e}rrez, M.~R. Orellana, A.~C. Robleda, V.~M.
  P{\'e}rez-Garc{\'\i}a, M.~Rosa, Dynamical properties of feedback signalling
  in b lymphopoiesis: A mathematical modelling approach, Journal of Theoretical
  Biology 522 (2021) 110685.

\bibitem{marciniak2019}
T.~Lorenzi, A.~Marciniak-Czochra, T.~Stiehl, A structured population model of
  clonal selection in acute leukemias with multiple maturation stages, Journal
  of mathematical biology 79~(5) (2019) 1587--1621.

\bibitem{clapp2015review}
G.~Clapp, D.~Levy, A review of mathematical models for leukemia and lymphoma,
  Drug Discovery Today: Disease Models 16 (2015) 1--6.

\bibitem{ducrot2007model}
A.~Ducrot, V.~Volpert, On a model of leukemia development with a spatial cell
  distribution, Mathematical Modelling of Natural Phenomena 2~(3) (2007)
  101--120.

\bibitem{moricke2005prognostic}
A.~M{\"o}ricke, M.~Zimmermann, A.~Reiter, H.~Gadner, E.~Odenwald, J.~Harbott,
  W.-D. Ludwig, H.~Riehm, M.~Schrappe, Prognostic impact of age in children and
  adolescents with acute lymphoblastic leukemia: data from the trials all-bfm
  86, 90, and 95, Klinische P{\"a}diatrie 217~(06) (2005) 310--320.

\bibitem{dai2021clinical}
Q.~Dai, G.~Zhang, H.~Yang, Y.~Wang, L.~Ye, L.~Peng, R.~Shi, S.~Guo, J.~He,
  Y.~Jiang, Clinical features and outcome of pediatric acute lymphoblastic
  leukemia with low peripheral blood blast cell count at diagnosis, Medicine
  100~(4) (2021).

\bibitem{bullhallmarks2022}
J.~A. Bull, H.~M. Byrne, The hallmarks of mathematical oncology, Proceedings of
  the IEEE (2022).

\bibitem{pui2006treatment}
C.-H. Pui, W.~E. Evans, Treatment of acute lymphoblastic leukemia, New England
  Journal of Medicine 354~(2) (2006) 166--178.

\bibitem{ronghe2001remission}
M.~Ronghe, G.~Burke, S.~Lowis, E.~Estlin, Remission induction therapy for
  childhood acute lymphoblastic leukaemia: clinical and cellular pharmacology
  of vincristine, corticosteroids, l-asparaginase and anthracyclines, Cancer
  treatment reviews 27~(6) (2001) 327--337.

\bibitem{jayachandran2014optimal}
D.~Jayachandran, A.~E. Rundell, R.~E. Hannemann, T.~A. Vik, D.~Ramkrishna,
  Optimal chemotherapy for leukemia: a model-based strategy for individualized
  treatment, PloS one 9~(10) (2014) e109623.

\bibitem{mouser2014model}
C.~L. Mouser, E.~S. Antoniou, J.~Tadros, E.~K. Vassiliou, A model of
  hematopoietic stem cell proliferation under the influence of a
  chemotherapeutic agent in combination with a hematopoietic inducing agent,
  Theoretical Biology and Medical Modelling 11~(1) (2014) 1--14.

\bibitem{nave2022new}
O.~Nave, A new protocol applied to cancer treatment-mathematical model-singular
  perturbed vector field algorithm (2022).

\bibitem{rubinow1976mathematical}
S.~Rubinow, J.~Lebowitz, A mathematical model of the chemotherapeutic treatment
  of acute myeloblastic leukemia, Biophysical journal 16~(11) (1976)
  1257--1271.

\bibitem{pefani2014chemotherapy}
E.~Pefani, N.~Panoskaltsis, A.~Mantalaris, M.~C. Georgiadis, E.~N.
  Pistikopoulos, Chemotherapy drug scheduling for the induction treatment of
  patients with acute myeloid leukemia, IEEE Transactions on Biomedical
  Engineering 61~(7) (2014) 2049--2056.

\bibitem{cartcells}
{\'A}.~Mart{\'\i}nez-Rubio, S.~Chuli{\'a}n, C.~Bl{\'a}zquez~Go{\~n}i,
  M.~Ram{\'\i}rez~Orellana, A.~P{\'e}rez~Mart{\'\i}nez, A.~Navarro-Zapata,
  C.~Ferreras, V.~M. P{\'e}rez-Garc{\'\i}a, M.~Rosa, A mathematical description
  of the bone marrow dynamics during {CAR T}-cell therapy in {B}-cell childhood
  acute lymphoblastic leukemia, International Journal of Molecular Sciences
  22~(12) (2021) 6371.

\bibitem{perez2021car}
V.~M. P{\'e}rez-Garc{\'\i}a, O.~Le{\'o}n-Triana, M.~Rosa,
  A.~P{\'e}rez-Mart{\'\i}nez, {CAR T} cells for {T}-cell leukemias: Insights
  from mathematical models, Communications in Nonlinear Science and Numerical
  Simulation 96 (2021) 105684.

\bibitem{leon2021car}
O.~Le{\'o}n-Triana, S.~Sabir, G.~F. Calvo, J.~Belmonte-Beitia, S.~Chuli{\'a}n,
  {\'A}.~Mart{\'\i}nez-Rubio, M.~Rosa, A.~P{\'e}rez-Mart{\'\i}nez,
  M.~Ramirez-Orellana, V.~M. P{\'e}rez-Garc{\'\i}a, {CAR T} cell therapy in
  {B}-cell acute lymphoblastic leukaemia: Insights from mathematical models,
  Communications in Nonlinear Science and Numerical Simulation 94 (2021)
  105570.

\bibitem{kimmel2021roles}
G.~J. Kimmel, F.~L. Locke, P.~M. Altrock, The roles of {T} cell competition and
  stochastic extinction events in chimeric antigen receptor {T} cell therapy,
  Proceedings of the Royal Society B 288~(1947) (2021) 20210229.

\bibitem{hallmarks2011}
D.~Hanahan, R.~A. Weinberg, Hallmarks of cancer: the next generation, Cell
  144~(5) (2011) 646--674.

\bibitem{karon1966vincris}
M.~Karon, E.~J. Freireich, E.~Frei~III, R.~Taylor, I.~J. Wolman, I.~Djerassi,
  S.~L. Lee, A.~Sawitsky, J.~Hananian, O.~Selawry, et~al., The role of
  vincristine in the treatment of childhood acute leukemia, Clinical
  Pharmacology \& Therapeutics 7~(3) (1966) 332--339.

\bibitem{campana2010minimal}
D.~Campana, Minimal residual disease in acute lymphoblastic leukemia,
  Hematology 2010, the American Society of Hematology Education Program Book
  2010~(1) (2010) 7--12.

\bibitem{drugbankPrednisona}
D.~Czock, F.~Keller, F.~M. Rasche, U.~H{\"a}ussler, Pharmacokinetics and
  pharmacodynamics of systemically administered glucocorticoids, Clinical
  pharmacokinetics 44~(1) (2005) 61--98.

\bibitem{drugbankVincri}
M.~Qweider, J.~M. Gilsbach, V.~Rohde, Inadvertent intrathecal vincristine
  administration: a neurosurgical emergency: case report, Journal of
  Neurosurgery: Spine 6~(3) (2007) 280--283.

\bibitem{drugbankDauno}
F.~M. Balis, J.~S. Holcenberg, W.~A. Bleyer, Clinical pharmacokinetics of
  commonly used anticancer drugs, Clinical pharmacokinetics 8~(3) (1983)
  202--232.

\bibitem{drugbankAspa}
B.~Asselin, C.~Rizzari, Asparaginase pharmacokinetics and implications of
  therapeutic drug monitoring, Leukemia \& lymphoma 56~(8) (2015) 2273--2280.

\bibitem{armstrong1987applications}
G.~M. Armstrong, C.~P. Midgley, Applications: The exponential-decay law applied
  to medical dosages., Mathematics Teacher 80~(2) (1987) 110--13.

\bibitem{haycock1978geometric}
G.~B. Haycock, G.~J. Schwartz, D.~H. Wisotsky, Geometric method for measuring
  body surface area: a height-weight formula validated in infants, children,
  and adults, The Journal of pediatrics 93~(1) (1978) 62--66.

\bibitem{szczepanski2007mdr}
T.~Szczepa{\'n}ski, Why and how to quantify minimal residual disease in acute
  lymphoblastic leukemia?, Leukemia 21~(4) (2007) 622--626.

\bibitem{amin2005having}
H.~Amin, Y.~Yang, Y.~Shen, E.~Estey, F.~Giles, S.~Pierce, H.~Kantarjian,
  S.~O'Brien, I.~Jilani, M.~Albitar, Having a higher blast percentage in
  circulation than bone marrow: clinical implications in myelodysplastic
  syndrome and acute lymphoid and myeloid leukemias, Leukemia 19~(9) (2005)
  1567--1572.

\bibitem{bianconi2013estimation}
E.~Bianconi, A.~Piovesan, F.~Facchin, A.~Beraudi, R.~Casadei, F.~Frabetti,
  L.~Vitale, M.~C. Pelleri, S.~Tassani, F.~Piva, et~al., An estimation of the
  number of cells in the human body, Annals of human biology 40~(6) (2013)
  463--471.

\bibitem{choi2014reference}
J.~Choi, S.~J. Lee, Y.~A. Lee, H.~G. Maeng, J.~K. Lee, Y.~W. Kang, Reference
  values for peripheral blood lymphocyte subsets in a healthy korean
  population, Immune network 14~(6) (2014) 289--295.

\bibitem{raes2006reference}
A.~Raes, S.~Van~Aken, M.~Craen, R.~Donckerwolcke, J.~V. Walle, A reference
  frame for blood volume in children and adolescents, BMC pediatrics 6~(1)
  (2006) 1--8.

\bibitem{linderkamp1977estimation}
O.~Linderkamp, H.~Versmold, K.~Riegel, K.~Betke, Estimation and prediction of
  blood volume in infants and children, European journal of pediatrics 125~(4)
  (1977) 227--234.

\bibitem{acharya1984development}
R.~Acharya, M.~K. Sundareshan, Development of optimal drug administration
  strategies for cancer-chemotheraphy in the framework of systems theory,
  International journal of bio-medical computing 15~(2) (1984) 139--150.

\bibitem{kay2013improving}
K.~Kay, I.~M. Hastings, Improving pharmacokinetic-pharmacodynamic modeling to
  investigate anti-infective chemotherapy with application to the current
  generation of antimalarial drugs, PLoS computational biology 9~(7) (2013)
  e1003151.

\end{thebibliography}


\begin{thebibliography}{[99]}
\bibitem{CRLS} D.P. Carlisle, \textsl{Packages in the `graphics' bundle},
available from CTAN as \verb=grfguide.tex= and  \verb=grfguide.ps=.

\bibitem{GMS} Michel Goossens, F. Mittelbach, and A.
Samarin, \textsl{The \LaTeX\ companion}, Addison-Wesley Co., Reading,
MA, 1994.

\bibitem{TB} D.E. Knuth, \textsl{The \TeX book},
Addison-Wesley, Reading, MA, 1984.

\bibitem{LM} L. Lamport, \textsl{\LaTeX: A document preparation
system}, 2nd revised ed., Addison-Wesley, Reading, MA, 1994.

\bibitem{RECK} K. Reckdahl, \textsl{Using EPS Graphics in \LaTeX2${}_\epsilon$ Documents},
available as \verb=epslatex.ps= from the \\
\verb=ftp://ftp.tex.ac.uk/tex-archive/info/= directory (or other CTAN sites).

\bibitem{RS} R. S\'eroul, \textsl{Le petit livre de \TeX},
Masson, Paris, 1996.

\bibitem{Joy} M.D. Spivak, \textsl{The joy of \TeX},
2nd revised ed., Amer. Math. Soc., Providence, RI, 1990.

\end{thebibliography}

\begin{comment}
HAY QUE PONERLO EN ESTE FORMATO. A PARTIR DE AHORA

\end{comment}

\end{document}

\maketitle

\footnotesize{
\noindent $^{1}$ \quad Department of Mathematics, Universidad de C\'{a}diz, Puerto Real, C\'{a}diz, Spain\\
$^{2}$ \quad Biomedical Research and Innovation Institute of C\'adiz (INiBICA), Hospital Universitario Puerta del Mar, C\'{a}diz, Spain\\
$^{3}$ \quad Department of Pediatric Hematology and Oncology, Hospital de Jerez C\'adiz, Spain\\
$*$ ana.nino@uca.es
}
\section*{Theoretical results}
\subsection*{Treatment Model}
\begin{equation}
\label{eq:ModeloTratamiento}
\begin{align}
\dfrac{dC_1}{dt}&= c_0 + s\rho_1C_1-\alpha_1C_1-\mu\rho_1C_1,\\
\dfrac{dC_2}{dt}&= s\rho_2C_2+\alpha_1C_1-\alpha_2C_2-\mu\rho_2C_2,\\
\dfrac{dC_3}{dt}&= \alpha_2C_2-\alpha_3C_3,\\
\dfrac{d L}{d t}&= s_{L} \rho_{L} L\left(1-\dfrac{L}{L_{\max }}\right)-\gamma_{L} L-\mu \rho_L L,\\
%las ecuaciones del tto
\mu&= \sum_{j\in J}\delta_j \left( \mu_j + Q_j\right),\\
Q_j&=\begin{cases} 
      q_j& t \in \mathcal{D}_j,\\
      &\\
      0 & t \notin\mathcal{D}_j,\\
   \end{cases},\\
 \dfrac{d\mu_j}{dt}&=
      -\lambda_j \mu_j(t).
\end{align}
\end{equation}
Note that to study Model~\eqref{eq:ModeloTratamiento}, treatments and dose assumptions have been simplified. In fact, Model~\eqref{eq:ModeloTratamiento} begins to stabilize from the end of the treatment due to the type of function $\mu$. In this way, two situations could be considered: \textbf{a)} study of the steady states from day 38 with initial data from this day \textbf{b)} Analyze the problem for a constant value $\mu=\Tilde{\mu}$. In this case, $\Tilde{\mu}$ will be a mean of all $\mu$ values.

\begin{enumerate}[label=\textbf{\alph*)}]
    \item If we consider the situation from day 38, leukemic cells are 0 and $\mu=0$. Therefore, it would be Model presented in healthy bone marrow and it will have the same stability values.

    \item To analyze the problem with $\mu=\Tilde{\mu}$, we consider $\Tilde{s}=s-\Tilde{\mu}$. Thus, we will have a system of equations analogous to Leukemic Model.
\end{enumerate}

\noindent Neither case provides information for the numerical study of the model.\\

%\noindent Generally, if theoretical results are proved for Model ~\eqref{eq:ModeloTratamiento}, consequently these will be proved for Model ~\eqref{eq:ModelosClon}, due to the fact that Model ~\eqref{eq:ModelosClon} is equal to the particular case to Model ~\eqref{eq:ModeloTratamiento} for $\mu=0$.

\subsubsection*{Existence, boundedness and positivity of solutions}

\textbf{Existence.}
Let us consider $C_1,C_2,C_3,L>0$ and $C_1(t_0)=C_1^0$, $C_2(t_0)=C_2^0$, $C_3(t_0)=C_3^0$, $L(t_0)=L^0$. The initial value problem for Model \eqref{eq:ModeloTratamiento} has a unique local-in-time solution for each $t \in \left[t_0-\epsilon, t_0+\epsilon \right]$ for some $\epsilon>0$.

\noindent \textit{Proof}.
Taking into account the continuity of the functions in Leukemic Model, along with the continuity of function $\mu$, it is known that at least there is a solution for Model~\eqref{eq:ModeloTratamiento} for each $(C_1^0,C_2^0,C_3^0,L^0)\in \mathbb{R}^4$ as initial values. \\ 

\noindent Moreover, boundedness of the partial derivatives of $C_i$ for $i=1,2,3$ and $L$, implies that they satisfy Lipschitz conditions, therefore, by the Picard-Lindelöf theorem, solutions of Model ~\eqref{eq:ModeloTratamiento} with initial values $C_1(t_0)=C_1^0$, $C_2(t_0)=C_2^0$, $C_3(t_0)=C_3^0$, $L(t_0)=L^0$ are unique.\\

\noindent \textbf{Positivity.}
The solutions of Model \eqref{eq:ModeloTratamiento} with initial value $\left(C_1^0,C_2^0,C_3^0,L^0\right)$ are positive.\\

\noindent\textit{Proof.}
Firstly, taking into account other proofs, $\left[26\right]$, it is known that solutions $C_i$ are positive. Considering that 
\begin{subequations}
    \begin{align}
    \dfrac{dC_i}{dt}&>-\alpha_i C_i-\mu \rho_i C_i.
    \end{align}
If both sides of the equation are integrated from $t_0$ to $t$, we obtain
    \begin{align}
    C_i(t)&>C_i^0\exp\left(\left(-\alpha_i -\mu \rho_i\right)t\right),
    \end{align}
so $C_i>0$. Therefore, let us study the leukemic compartment
    \begin{align}
\dfrac{dL}{dt}&=L(s_L\rho_L(1-\dfrac{L}{L_{max}})-\gamma_L-\mu \rho_i).
    \end{align}

Due to the fact that 
    \begin{align}
\dfrac{dL}{dt}&>\dfrac{-s_L \rho_L L^2}{L_{\max}}-\gamma_L L-\mu\rho_L L,
  \end{align}

if that expression is integrated,
\begin{align}
L(t)&>\dfrac{L_{\max}L_0(\gamma_L+\mu\rho_L)}{\rho_L s_L L_0 (e^{t (\gamma_L +\mu\rho_L)}-1)+ e^{t(\gamma_L +\mu \rho_L)}\gamma_L L_{\max}}.
 \end{align}
\end{subequations}

The numerator is positive and also the denominator since $t\left({\gamma_L+\mu \rho_L}\right)>0$, as a result, $e^{t\left({\gamma_L+\mu \rho_L}\right)}>1$.\\

\noindent \textbf{Boundedness.}
Solutions for Model~\eqref{eq:ModeloTratamiento} are bounded.
\textit{Proof.}
According to the previous study $\left[26\right]$,  it has just been proved for $C_i$ population that $C_i$, $i=1,2,3$ is bounded at the top, with the result that adding the term $-\mu\rho_iC_i$ which indicates a decrease, both $C_1$ and $C_2$ remain bounded. \\

\noindent Secondly, leukemic cells behavior $L$, has been defined by a logistic equation, therefore, it is bounded by the parameter $L_{\max}$.

\subsection*{Leukemia Model}
\label{appendix:theoricalSolLEUKEMIA}

\begin{equation}
\label{eq:ModelosClon}
\begin{align}
\dfrac{dC_1}{dt}&= c_0 + s\rho_1C_1-\alpha_1C_1,\\
\dfrac{dC_2}{dt}&= s\rho_2C_2+\alpha_1C_1-\alpha_2C_2,\\
\dfrac{dC_3}{dt}&= \alpha_2C_2-\alpha_3C_3,\\
\dfrac{d L}{d t} &= s_{L} \rho_{L} L\left(1-\dfrac{L}{L_{\max }}\right)-\gamma_{L} L,
\end{align}
\end{equation}

For simplicity, in the study of the Model \eqref{eq:ModelosClon}, the influx parameter $c_0$ is omitted without loss of generality. 
Let us study Model which describes Leukemic cells development, where any treatment is considered, Model ~\eqref{eq:ModelosClon}. 
Therefore, to obtain steady states for this model we assume non-existent flux in each stage, thus

\begin{subequations}
\begin{align}
0&=\dfrac{\rho_1C_1}{1+k\left(L+\sum_{i=1}^3{C_i}\right)}-\alpha_1C_1,\\
0&=\dfrac{\rho_2C_2}{1+k\left(L+\sum_{i=1}^3{C_i}\right)}+\alpha_1C_1-\alpha_2C_2,\\
0&=\alpha_2C_2-\alpha_3C_3,\\
0&=\dfrac{\rho_{L} L}{1+k\sum_{i=1}^3{C_i}} \left(1-\dfrac{L}{L_{\max }}\right)-\gamma_{L} L.
%s_{L}(t)& = &\dfrac{1}{1+k\left(\sum_{i=1}^{3} C_{i}\right)}\\
%&&\\
\end{align}
\end{subequations}

\vspace{2em}
\noindent Consequently, we have six solutions to the above-mentioned system: \\

\vspace{1ex}
\scalebox{0.9}{
\[ \begin{array}{ll}
\hspace{-3em} P_{L_1}=&\left(0,0,0,0\right) \\
&\\
\hspace{-3em} P_{L_2}=&\left(0,0,0, \dfrac{ L_{\max}(\rho_L -\gamma_L )}{\rho_L} \right)\\
&\\
\hspace{-3em} P_{L_3}=&\left(0, \dfrac{\alpha_3(\rho_2-\alpha_2)}{\alpha_2(\alpha_2+\alpha_3)k},\dfrac{\rho_2-\alpha_2}{(\alpha_2+\alpha_3)k},0 \right) \\
&\\
\hspace{-3em} P_{L_4}=&\left(  -\dfrac{\alpha_3(\alpha_2\rho_1-\alpha_1\rho_2)\omega}{ \alpha_1 k \phi \psi}    ,   \dfrac{-\alpha_1 \alpha_3 \rho_1 +\alpha_3\rho_1^2}{ k \phi}   ,\dfrac{-\alpha_1 \alpha_2 \rho_1 +\alpha_2\rho_1^2}{ k \phi},0\right) \\
&\\
\hspace{-3em} P_{L_5}=&\left(-\dfrac{\alpha_3(\alpha_2\rho_1-\alpha_1\rho_2)\omega}{ \alpha_1 k \phi \psi},-\dfrac{\alpha_3\rho_1\omega}{k \phi\psi},-\dfrac{\alpha_2\rho_1\omega}{k\phi\psi},\dfrac{L_{\max}\rho_1(\gamma_L -\alpha_1)}{\alpha_1\psi} \right)\\
&\\
\hspace{-3em} P_{L_6}=&\left(0,\alpha_2\alpha_3\dfrac{-\gamma_LkL_{\max}+\rho_L+kL_{\max}\rho_L-\dfrac{\rho_2\rho_L}{\alpha_2}}{\alpha_2\left(\alpha_2+\alpha_3\right)k\psi}  ,
\dfrac{-\rho_2\rho_L+\alpha_2\left(-\psi+kL_{\max}\rho_L\right) }{(\alpha_2+\alpha_3 )k \psi},\dfrac{\gamma_L L_{\max}\rho_2-\alpha_2 L_{\max}\rho_L}{\alpha_2 \psi} \right) \\
 \end{array}\] 
}

\vspace{2em}
Being $\phi=\alpha_1\alpha_2\rho_1+\alpha_1 \alpha_3 \rho_1 + \alpha_2\alpha_3 \rho_1 -\alpha_1 \alpha_3 \rho_2$, $ \psi=\gamma_L k L_{\max}-\rho_L $ and $\omega=\alpha_1 \gamma_L k L_{\max}-\alpha_1 \rho_L-\alpha_1 k L_{\max}\rho_L+\rho_1\rho_L $\\

Firstly, values obtained as solutions must be positive so that they are biologically well defined.\\

The Jacobian matrix related to the system at any point $\left(C_1,C_2,C_3,L\right)$ is\\

\begin{equation}

\left(
\begin{array}{cccc}
-\alpha_1-\dfrac{C_1 k \rho_1}{\mathcal{A}^2}+\dfrac{ \rho_1}{\mathcal{A}}& -\dfrac{C_1 k \rho_1}{\mathcal{A}^2} & -\dfrac{C_1 k \rho_1}{\mathcal{A}^2}&-\dfrac{C_1 k \rho_1}{\mathcal{A}^2} \\
&&&\\
 -\alpha_1-\dfrac{C_2 k \rho_2}{\mathcal{A}^2}& -\alpha_2-\dfrac{C_2 k \rho_2}{\mathcal{A}^2}+\dfrac{ \rho_2}{\mathcal{A}}  & -\dfrac{C_2 k \rho_2}{\mathcal{A}^2}& -\dfrac{C_2 k \rho_2}{\mathcal{A}^2} \\
&&&\\
 0&  \alpha_2&-\alpha_3 &0 \\
&&&\\
 -\dfrac{Lk\mathcal{C}}{\mathcal{B}}&   -\dfrac{Lk\mathcal{C}}{\mathcal{B}}& -\dfrac{Lk\mathcal{C}}{\mathcal{B}} & -\gamma_L + \dfrac{\mathcal{C}}{\mathcal{B}}-\dfrac{\left(L\rho_L\right)\rho_L}{\mathcal{B}L_{\max}}\\
&&&\\
\end{array}
\right)
\label{JacobianaL}
\end{equation}

\vspace{3em}

\noindent Being, $\mathcal{A}=1+\left(1+C_1+C_2+C_3+L\right)k$, \hspace{1em}
$\mathcal{B}=1+\left(1+C_1+C_2+C_3\right)k$ \hspace{1em} \\

\vspace{1 ex}
\noindent and $\mathcal{C}=\dfrac{-\left(1-\frac{L}{L_{\max}}\right)\rho_L}{\mathcal{B}}$.\\

Substituting each $P_{L_i}$, $i=1,\dots,6$ in (\ref{JacobianaL}) we get its own  eigenvalues. Thus, we can propose stability conditions thanks to Hartmann-Grobmann Theorem.\\

\subsection*{Leukemia Model Numerical Results}
\label{appendix:numericalLEUKEMIA}
According to theoretical results along with data in Table Parameters Value, it is possible to find steady states for Model ~\eqref{eq:ModelosClon}. Therefore, thanks to Hartman-Grobman Theorem, steady states can be classified by stability in each case.  \\

Assuming that both cases ($\rho_L=\rho_1$ and $\rho_L=\rho_2$) are analogous, we present the study for case I ($\rho_L=\rho_1$). We consider the system:
\begin{subequations}
\label{eq:ClonNum}
\begin{align}
f_1\left(C_1,C_2,C_3,L\right)&= c_0 + s\rho_1C_1-\alpha_1C_1,\\[5pt]
f_2\left(C_1,C_2,C_3,L\right)&= s\rho_2C_2+\alpha_1C_1-\alpha_2C_2,\\[5pt]
f_3\left(C_1,C_2,C_3,L\right)&= \alpha_2C_2-\alpha_3C_3,\\[2pt]
f_4\left(C_1,C_2,C_3,L\right)&= s_{L} \rho_{L} L\left(1-\dfrac{L}{L_{\max }}\right)-\gamma_{L} L.
\end{align}
\end{subequations}

Firstly, we obtain the steady states from solving $f_i=0$, for $i=1,\dots,4$ in system ~\eqref{eq:ClonNum}. Thus, there are six steady states presented in Table \ref{table:eigenvalues}. Moreover, the Jacobian matrix for the system \eqref{eq:ClonNum} is the same as the Jacobian matrix \ref{JacobianaL} where we can substitute the parameter values from Table Parameters Value. To study the stability of each steady state, just replace the value of each $P_{L_i}$ above in Jacobian matrix \ref{JacobianaL}. Then, we calculate the eigenvalues corresponding to these arrays. We get the values presented in the Table ~\ref{table:eigenvalues} and consequently, the stability obtained for each point according to Hartman-Grobman Theorem.

\begin{flalign*}
P_{L_1}=&\left(-1.13654 \times 10^{11}\right.,&-5.97781  \times 10^{11}&,&-2.98891  \times 10^{11}&,&\left.1.04156 \times 10^{12} \right) & \\
P_{L_2}=&\left(-2.0867 \times 10^{8}\right.,&-6.79766 \times 10^{11}&,&-3.39883 \times 10^{11}&,& \left.1.04196\times10^{12} \right) &\\
P_{L_3}=&\left(-1.94642 \times 10^{8}\right. ,&1.45273\times10^{10}&,&7.26363\times10^{9}&,&\left.0 \right) &\\
P_{L_4}=&\left(-1.90323\times10^7\right. ,&1.00459\times10^{7}&,&5.02296\times10^6&,&\left.0\right) &\\
P_{L_5}=&\left(6.20596\times10^{7}\right.,&7.47794\times10^{7}&,&3.73897\times10^{7}&,&\left.9.99577\times10^{11}\right) &\\
P_{L_6}=&\left(3.77184\times10^{9}\right.,&1.87656\times10^{10}&,&9.38282\times10^{9}&,&\left.0 \right) &\\
\end{flalign*}

\begin{comment}

\begin{enumerate}[label= $P_{L_{\arabic*}}$]\itemsep=3em
    \item $\left(-1.13654 \times 10^{11},-5.97781  \times 10^{11},-2.98891  \times 10^{11},1.04156 \times 10^{12}\right)$ 
    \item 
    $\left(-2.0867 \times 10^{8},-6.79766 \times 10^{11},-3.39883 \times 10^{11}, 1.04196\times10^{12}\right)$ 
    \item 
    $\left(-1.94642 \times 10^{8},1.45273\times10^{10},7.26363\times10^{9},0\right)$ 
    \item 
    $\left(-1.90323\times10^7,1.00459\times10^{7},5.02296\times10^6,0 \right)$ 
     \item 
    $\left(6.20596\times10^{7},7.47794\times10^{7},3.73897\times10^{7},9.99577\times10^{11}\right)$ 
    \item 
    $\left(3.77184\times10^{9},1.87656\times10^{10},9.38282\times10^{9},0\right)$
   
\end{enumerate}

\end{comment}
\vspace{3em}

\begin{table}[H] 

%%% \tablesize{} %% You can specify the fontsize here, e.g., \tablesize{\footnotesize}. If commented out \small will be used.
\begin{adjustbox}{max width=1\textwidth,center}
\begin{tabular}{ccc}
\toprule
\textbf{Steady state}	& \textbf{Eigenvalues}	& Stability (H-G) \\
\midrule
$P_{L_{1}}$& $\left(2.19456,-0.381211,-0.045143,0.00691681\right)$ & Unstable\\
&&\\
$P_{L_{2}}$&$\left(3.17661,-0.414698,0.0478288,0.0068542\right)$ &Unstable\\
&&\\
$P_{L_{3}}$&$\left(-0.227986,0.219088,-0.12632,0.054041\right)$ &Unstable\\
&&\\
$P_{L_{4}}$&$\left(0.693134,0.52798,0.316404,-0.287825\right)$ &Unstable\\
&&\\
$P_{L_{5}}$&$\left(-0.680989,-0.288,-0.161129,-0.13943\right)$ &Stable\\
&&\\
$P_{L_{6}}$&$\left(-0.254221,0.165061,-0.0672105+0.0249579i, -0.0672105-0.0249579 i \right)$ &Unstable\\
\bottomrule
\end{tabular}
\end{adjustbox}
\caption{\textbf{Eigenvalues related to each steady state from system \eqref{eq:ClonNum}}. Values which are obtained from substituting each steady state in the Jacobian matrix associated to the Model \eqref{eq:ModelosClon}. }
\label{table:eigenvalues}
\end{table}

\vspace{2em}

As it is shown in Table ~\ref{table:eigenvalues}, there are five steady states which are unstable. Solutions $P_{L_1}, P_{L_2}, P_{L_3}$, and $P_{L_4}$ display negative values for some compartments, therefore, they are meaningless. There is a steady state, $P_{L_6}$, which represents the situation where there is no leukemia. This point is also unstable since it is related to the healthy bone marrow compartments, hence, if a leukemic cell appears, values change considerably. 
As a result, there is an unique stable steady state, $P_{L_{5}}$, which corresponds to the case where leukemic cells $L$ invade the bone marrow.

\clearpage

\section*{Parameters sensitivity}
\label{appendix:sensitivity}
\begin{comment}
    To study parameters influence in the Model~\eqref{eq:ModeloTratamiento}, a sensitivity analysis of values $\delta_j$ has been implemented. As expected, $\delta_P$ related to prednisone influence is the most influential parameter given that it is administered every days. 

\end{comment}

\label{appendix:param}

   \begin{figure}[h]
  \centering
  
     \begin{subfigure}[b]{0.4\textwidth}

          \caption{\\}
          \includegraphics[width=1.2\textwidth]{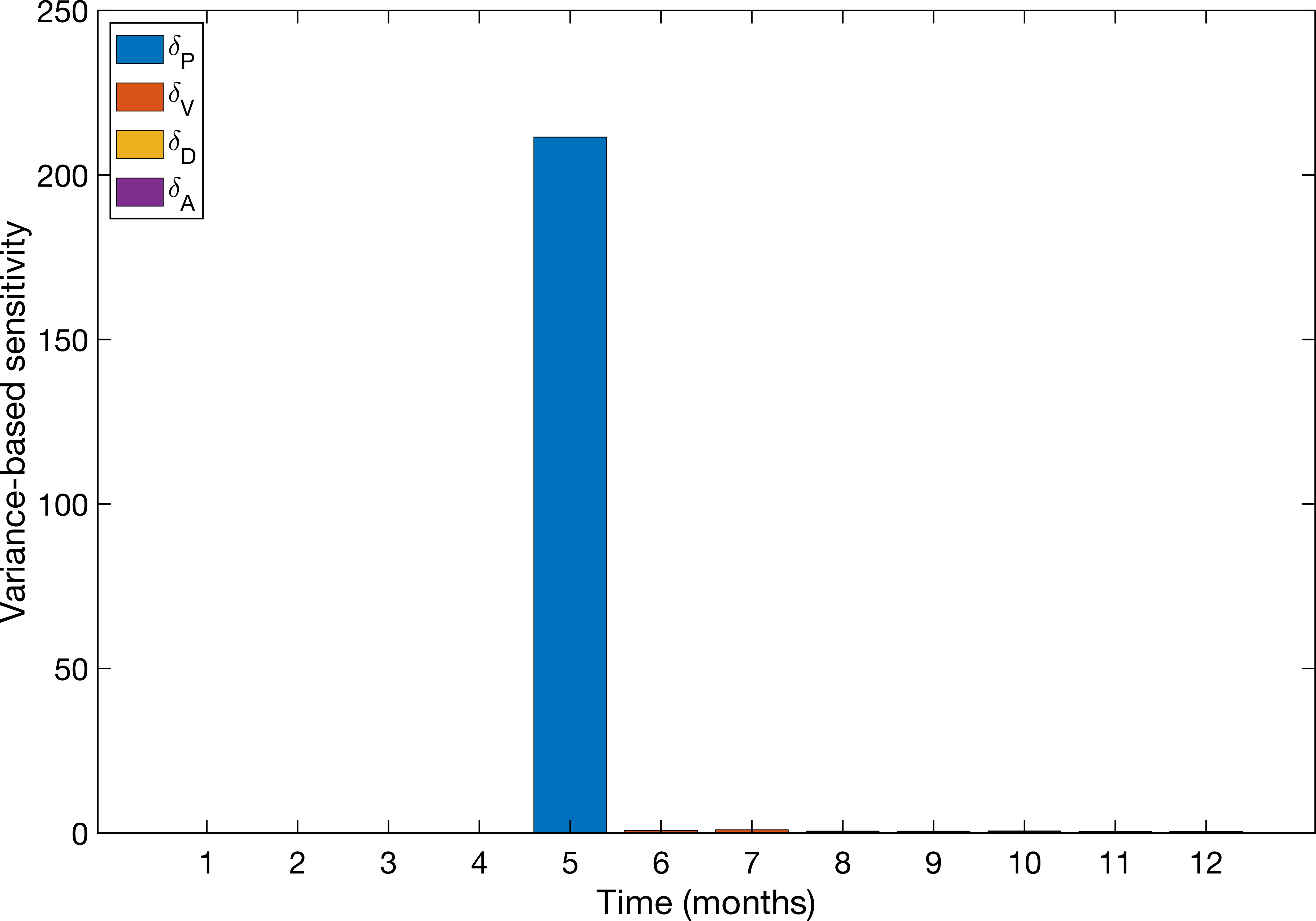}
          
        \end{subfigure}
      \begin{subfigure}[b]{0.4\textwidth}
      \caption{\\}
      \includegraphics[width=1.2\textwidth]{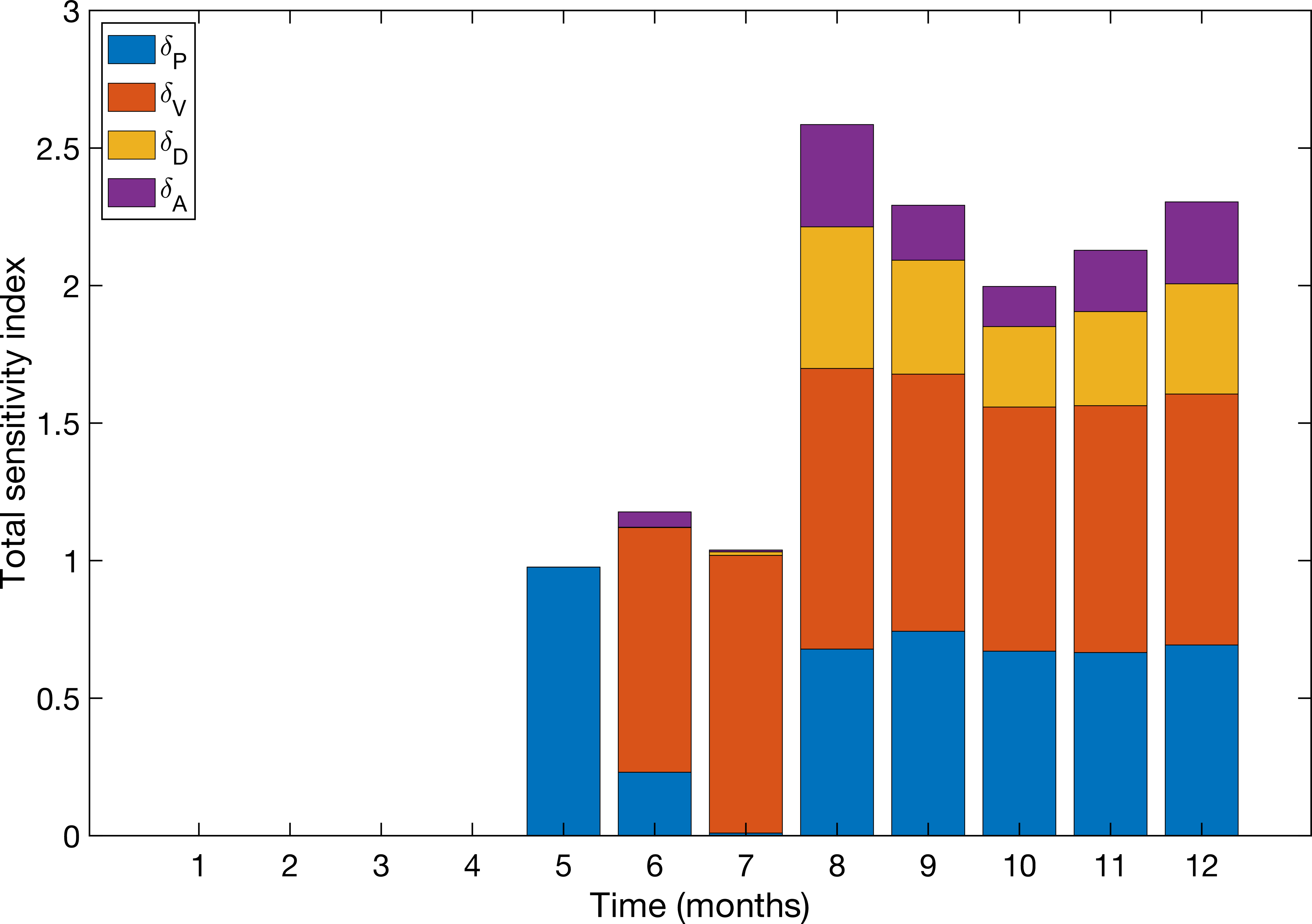}
      
    \end{subfigure}
    
     \begin{subfigure}[b]{0.4\textwidth}
          \caption{\\}
          
          \includegraphics[width=1.2\textwidth]{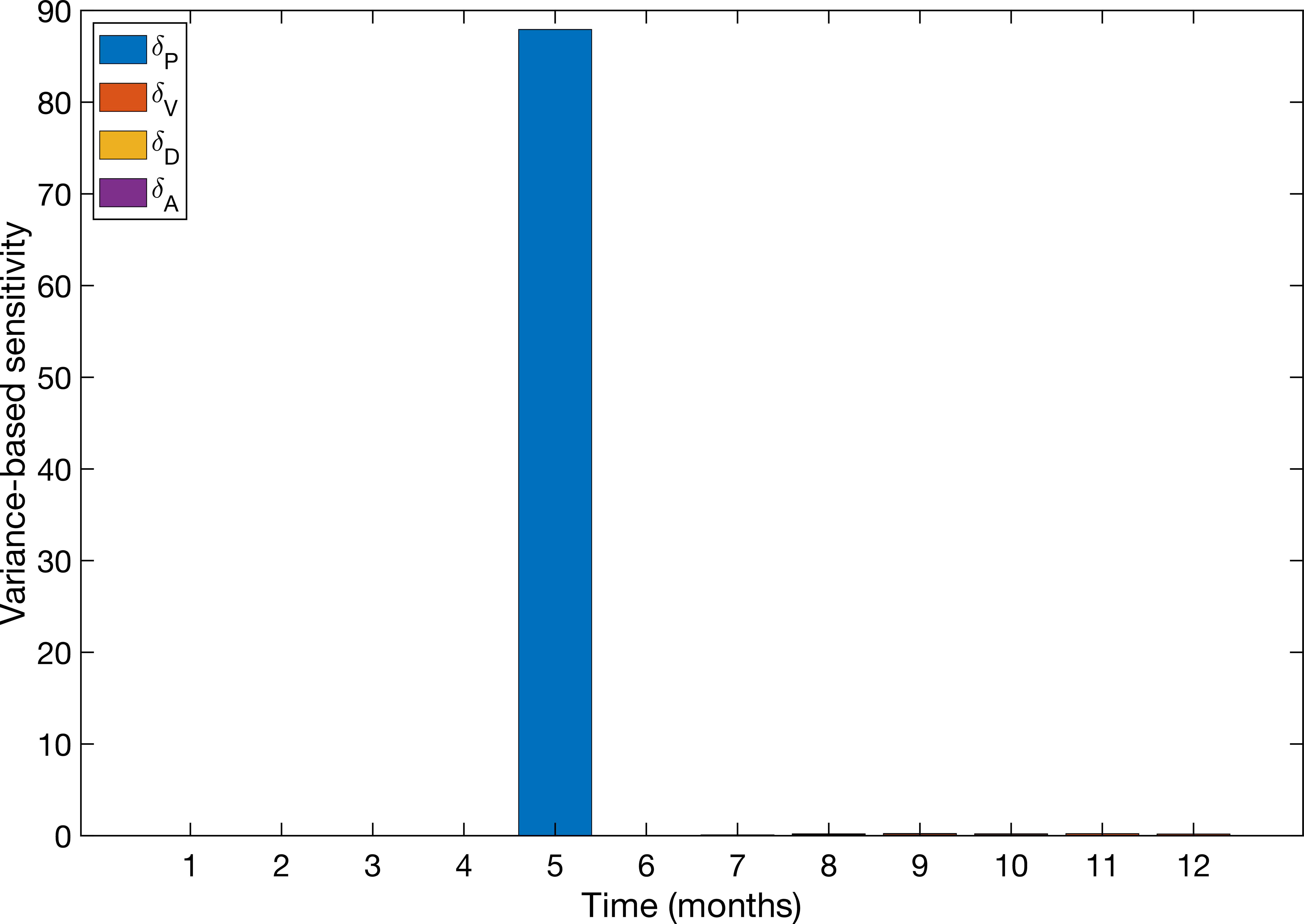}
          
        \end{subfigure}
      \begin{subfigure}[b]{0.4\textwidth}
      \caption{\\}
      \includegraphics[width=1.2\textwidth]{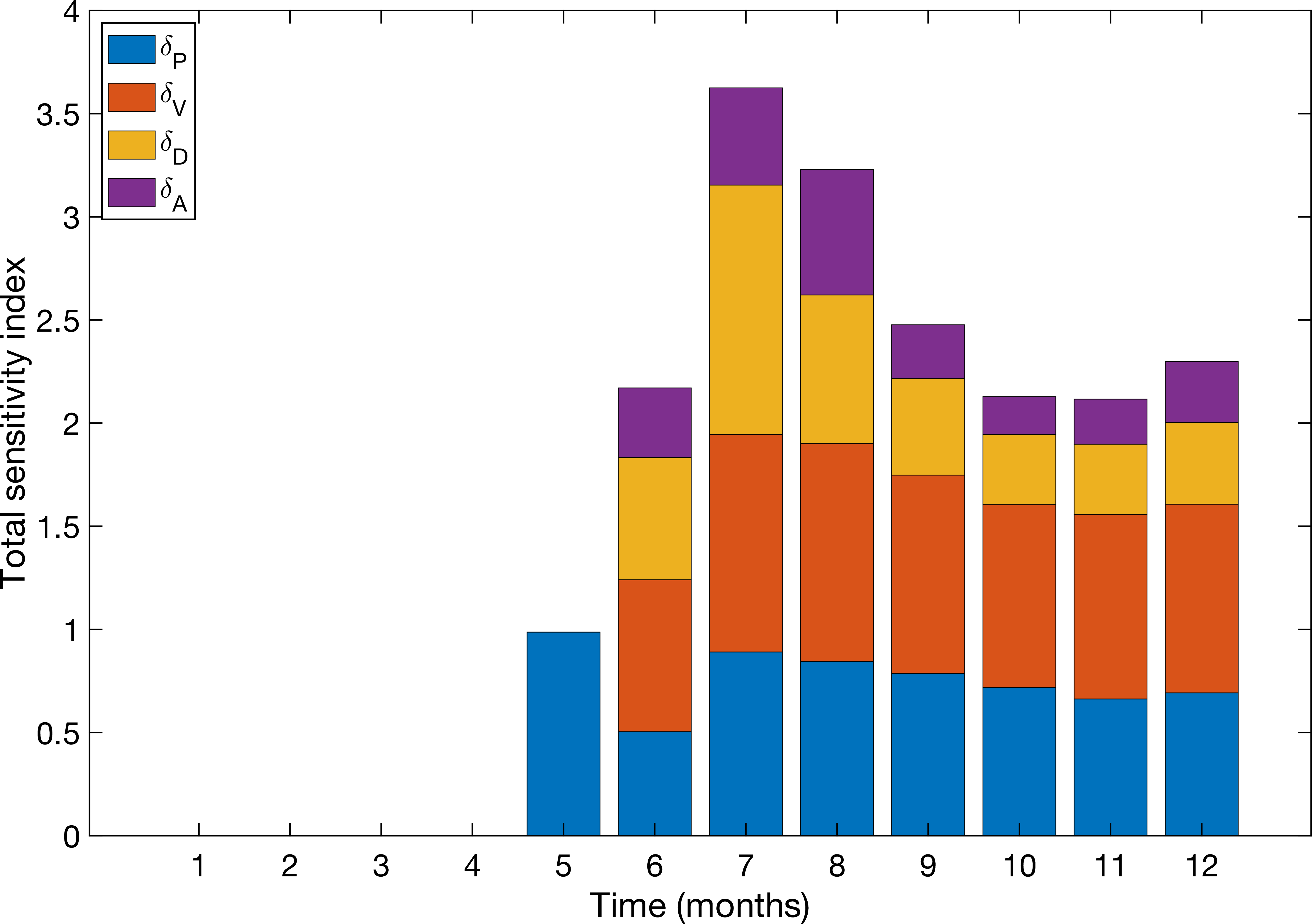}
      
    \end{subfigure}

     \caption{\textbf{Sobol's sensitivity analysis for $\delta_j$}. Contribution of drugs influence parameters, $\delta_j$, on the Model~\eqref{eq:ModeloTratamiento}. Parameters values in the defined intervals: $\delta_P \in \left[1/60, 1/6\right]$, $\delta_V \in \left[1/1.5, 10/1.5\right]$, $\delta_D \in \left[1/30, 1/3\right]$, $\delta_A \in \left[10^{-4}, 10^{-3}\right]$ \textbf{(a)} Variance-based sensitivity  on healthy cells\textbf{ (b)} Total influence on healthy cells \textbf{(c)} Variance-based sensitivity  on leukemic cells \textbf{(d)} Variance-based sensitivity  on leukemic cells }
     %Una era de las células leucémicas y otras de las sanas. Hablaría de: el más influyente es la prednisona y por tanto se ve en la primera que todo es par tal tal y claro, a partir del més cinco es porque es donde empieza el tratamiento, se ven influenciadas por prednisona, la otra y la otra. Y en los resultados en el mes 10 y 11 ya afectan todas por igual. por qué? porque no se dan o porque cuando se han dado se ha visto la influencia en el mes 10 y 11?
    \label{fig:sensitivity}
  \end{figure}

According to the fact that $\delta_j q_j$ have the same order of magnitude for all $j\in J$ along with Sobol's sensitivity analysis for those ranges of values, we approximate values of $\delta_j$.
If we simulate the treatment assuming one of them with its maximum possible value and the others equal to zero, there is not biological sense. After several simulations, we suggest that analysis of parameters begins from $\delta_j q_j = 2.5$ for prednisone and vincristine, fixing $\delta_D, \delta_A$ in those values due to their little influence on the total treatment.

\begin{figure}[H]
    \centering
          \includegraphics[width=1\textwidth]{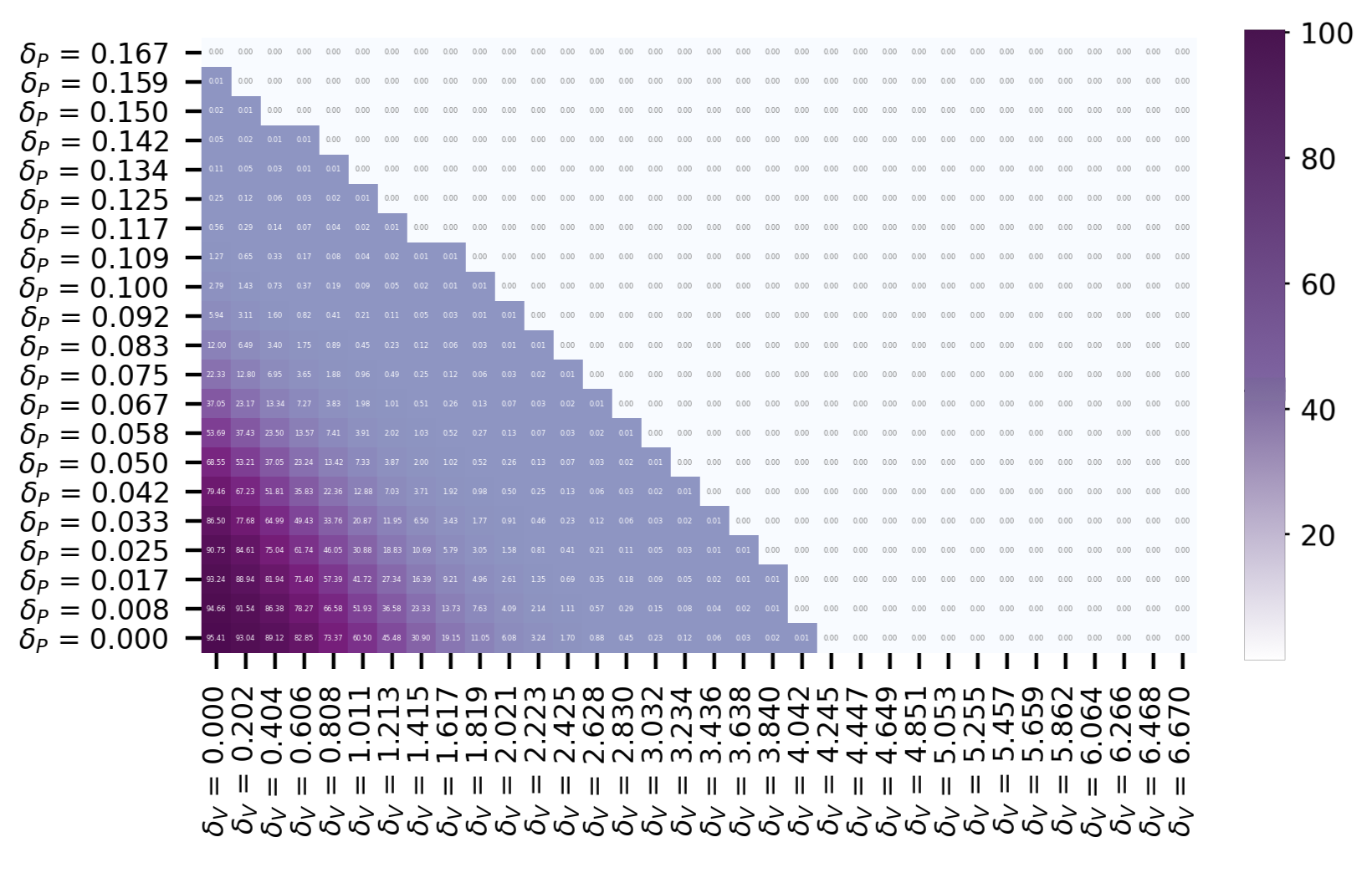}
     \caption{\textbf{Blasts percentage in a patient considering wider range for parameters variation.} We consider Prednisone values in the range $\left[1/60,1/6\right]$ and  the range Vincristine. For each values combination, Model \eqref{eq:ModeloTratamiento} is solved for parameters in Table Parameters Value with $\delta_D= 2.5/30 day/mg$, $\delta_A= 2.5/10^4 day/U$. The value in each box corresponds to the percentage of leukemic cells in day $+15$ depending on values $\delta_P$ and $\delta_V$. }
    \label{fig:FULL_heatmap}
  \end{figure}

\clearpage